\documentclass[11pt]{amsart}
\usepackage{amsmath}
\usepackage{amsxtra}
\usepackage{amscd}
\usepackage{amsthm}
\usepackage{amsfonts}
\usepackage{amssymb}
\usepackage{eucal}

\newcommand{\qtri}[4]{{\left[ 
\begin{array}{ccc}
 & \displaystyle #1 & \\ 
\displaystyle #2 & \displaystyle #3 & \displaystyle #4  
\end{array}
\right]
}}
\newcommand{\slth}{\widehat{\mathfrak{sl}}_2}
\newcommand{\slt}{\mathfrak{sl}_2}

\newcommand{\ch}{{\rm ch}}
\newcommand{\qbin}[2]{{\left[
\begin{matrix}{\displaystyle #1}\\
{\displaystyle #2}\end{matrix}
\right]
}}

\newcommand{\cF}{\mathcal{F}}
\newcommand{\cG}{\mathcal{G}}
\newcommand{\cH}{\mathcal{H}}
\newcommand{\cL}{\mathcal{L}}
\newcommand{\phih}{\hat{\phi}}
\newcommand{\tF}{\widetilde{\cF}}
\newcommand{\tG}{\widetilde{\cG}}
\newcommand{\tH}{\widetilde{\cH}}
\newcommand{\tN}{\widetilde{N}}
\newcommand{\nn}{\nonumber}
\newcommand{\bea}{\begin{eqnarray}}
\newcommand{\ena}{\end{eqnarray}}
\newcommand{\be}{\begin{eqnarray*}}
\newcommand{\en}{\end{eqnarray*}}
\def\bel{\begin{eqnarray}}
\def\enl{\end{eqnarray}}

\newcommand{\yc}{\check{y}}
\newcommand{\yt}{\tilde{y}}

\newcommand{\C}{{\mathbb C}}

\newcommand{\Z}{{\mathbb Z}} 
\newcommand{\al}{{\alpha}}

\newcommand{\la}{{\lambda}}

\def\r{\right}
\def\l{\left}
\def\P{\mathcal P}
\def\R{\mathcal R}
\numberwithin{equation}{section}

\numberwithin{equation}{section}
\newtheorem{thm}{Theorem}[section]
\newtheorem{prop}[thm]{Proposition}
\newtheorem{lem}[thm]{Lemma}
\newtheorem{cor}[thm]{Corollary}

\newtheorem{conj}[thm]{Conjecture}

\newcommand{\nc}{\newcommand}
\nc{\suth}{\widehat{\mathfrak{su}}(2)}
\nc{\gl}{\mathfrak{gl}}
\nc{\GL}{\mathfrak{GL}}
\nc{\g}{\mathfrak{g}}
\nc{\h}{\mathfrak{h}}
\nc{\na}{\mathfrak{n}}
\nc{\ve}{\varepsilon}
\nc{\Id}{{\mathop {\rm id}}}
\nc{\Tr}{{\mathop {\rm Tr}\,}}
\nc{\U}{{\mathop {\rm U}}}
\nc{\x}{{\bf x}}
\nc{\pa}{\partial}
\nc{\ld}{\ldots}
\nc{\cd}{\cdots}
\nc{\sm}{\sqrt{2m}\,}
\nc{\hk}{\hookrightarrow}
\nc{\A}{\mathfrak A}
\nc{\qb}[2]{\genfrac{(}{)}{0pt}{}{#1}{#2}_q}
\nc{\n}{\mathfrak{n}}
\nc{\un}{\mathfrak{u}}
\nc{\T}{\otimes}
\nc{\bv}{{\bf v}}
\nc{\bt}{{\bf t}}
\nc{\bu}{{\bf u}}
\nc{\bin}[2]{{\genfrac{[}{]}{0pt}{0}{#1}{#2}}_q}
\nc{\bininv}[2]{{\genfrac{[}{]}{0pt}{0}{#1}{#2}}_{q^{-1}}}
\nc{\fac}[1]{(#1)_q!}
\nc{\bl}{{\bf\lambda}}
\nc{\bs}{{\bf s}}
\nc{\bn}{{\bf n}}
\nc{\bm}{{\bf m}}

\renewcommand{\triangle}{\Delta}

\begin{document}

\title[(1,3)-filtration]
{A $\phi_{1,3}$-filtration of the Virasoro minimal 
series $M(p,p')$ \\
with $1< p'/p< 2$}
\author{B. Feigin, E. Feigin, M. Jimbo, T. Miwa and Y. Takeyama}
\address{BF: Landau institute for Theoretical Physics, 
Chernogolovka,
142432, Russia and \newline
Independent University of Moscow, Russia, Moscow, 119002, 
Bol'shoi Vlas'evski per., 11}
\email{feigin@mccme.ru}  
\address{EF: 
Tamm Theory Division, Lebedev Physics Institute, Russia, Moscow, 119991,
Leninski pr., 53 and\newline
Independent University of Moscow, Russia, Moscow, 119002, 
Bol'shoi Vlas'evski per., 11}
\email{evgfeig@mccme.ru}
\address{MJ: Graduate School of Mathematical Sciences, 
The University of Tokyo, Tokyo 153-8914, Japan}
\email{jimbomic@ms.u-tokyo.ac.jp}
\address{TM: Department of Mathematics, 
Graduate School of Science,
Kyoto University, Kyoto 606-8502,
Japan}\email{tetsuji@math.kyoto-u.ac.jp}
\address{YT: Institute of Mathematics, 
Graduate School of Pure and Applied Sciences, 
University of Tsukuba, Tsukuba, Ibaraki 305-8571,
Japan}\email{takeyama@math.tsukuba.ac.jp}

\dedicatory{Dedicated to Professor Heisuke Hironaka}

\date{\today}

\begin{abstract}
The filtration of the Virasoro minimal series representations
$M^{(p,p')}_{r,s}$
induced by the $(1,3)$-primary field $\phi_{1,3}(z)$
is studied. For $1< p'/p< 2$, a conjectural basis of $M^{(p,p')}_{r,s}$
compatible with the filtration is given by using monomial vectors in terms of
the Fourier coefficients of $\phi_{1,3}(z)$. In support of this conjecture,
we give two results. First, we establish the equality of the character
of the conjectural basis vectors with the character of the whole representation
space. Second, for the unitary series ($p'=p+1$), we establish for each $m$
the equality
between the character of the degree $m$ monomial basis and the character of the
degree $m$ component in the associated graded module
${\rm Gr}(M^{(p,p+1)}_{r,s})$ with respect to the filtration defined by
$\phi_{1,3}(z)$.
\end{abstract}
\maketitle

\section{Introduction}
In this paper we study a filtration 
on the Virasoro minimal modules by the 
$\phi_{1,3}$ primary field. 
We first state the problem in a general scheme.
Let $S=\{\phi_{\al}(z)\}_{\al\in A}$ 
be a set of vertex operators acting on a graded vector space $V$. 
In the actual setting, the representation space $V$ is a direct sum
$V=\oplus_{\beta\in B}V^{(\beta)}$ such that the index set $A$
is a subset of $B$, the grading is of the form
$V^{(\beta)}=\sum_{i\in\Delta^{(\beta)}+\Z_{\geq0}}V^{(\beta)}_i$,
and the action of the vertex operator
\be
\phi_\alpha(z):V^{(\beta)}\rightarrow V^{(\gamma)}
\en
is decomposed as $\sum_{n\in\Z+\Delta^{(\gamma)}-\Delta^{(\beta)}}
\phi^{(\gamma,\beta)}_{\alpha,-n}z^{n-\Delta^{(\alpha)}}$,
where $\phi^{(\gamma,\beta)}_{\alpha,-n}:V^{(\beta)}_i\rightarrow V^{(\gamma)}_{i+n}$.

Now for some fixed vector $v_0\in V^{(\alpha_0)}$,  
one can define a sequence of subspaces 
$E_0(V)\subset E_1(V)\subset E_2(V)\subset\cdots\subset V$ by 
setting 
\bea
& & E_m(V)
=\mathrm{span}
\{
\phi^{(\beta_0,\beta_1)}_{\al_1,-n_1}
\phi^{(\beta_1,\beta_2)}_{\al_2,-n_2}
\ldots \phi^{(\beta_{k-1},\beta_k)}_{\al_k,-n_k}v_0\mid\label{fil}\\
& & \phantom{quad}
\al_i\in A,\beta_j\in B,\beta_k=\alpha_0,\,
n_i\in\Z+\Delta^{(\beta_{i-1})}-\Delta^{(\beta_i)},\,k\le m
\}.\nonumber
\ena
In what follows we assume that 
the Fourier coefficients $\{\phi^{(\gamma,\beta)}_{\al,-n}\}$ 
generate the whole $V$ from $v_0$. 
In this case the above construction 
gives a filtration on $V$, which we refer to as 
the $S$-filtration.  
In our examples, 
$V$ is a representation of the Virasoro or 
the $\slth$ algebra, 
and $\phi_{\al}(z)$ 
are vertex operators from the corresponding conformal field theory.

There are two natural questions about $E_m(V)$:
\begin{itemize}
\item[(i)]
Find the $(q,v)$-character of $V^{(\alpha)}$:
$$
\ch_{q,v} V^{(\alpha)}=\sum_{m,n} 
q^n v^m \dim\bigl(V^{(\alpha)}_n\cap E_m(V))
/(V^{(\alpha)}_n\cap E_{m-1}(V)).
$$

\item[(ii)]
Find a monomial basis of $V$ 
which is compatible with $E_m(V)$. 
This means that one needs to construct a basis of 
$V$ of the form 
$$ 
\phi^{(\beta_0,\beta_1)}_{\al_1,-n_1}
\phi^{(\beta_1,\beta_2)}_{\al_2,-n_2}
\ldots \phi^{(\beta_{k-1},\beta_k)}_{\al_k,-n_k}v_0\quad(\beta_k=\alpha_0)
$$
with certain $\al_i, n_i,\beta_j$ such that the 
images of the basis vectors
with $k\le m$ form a basis of $E_m(V)$. 
\end{itemize}

In the case when $V$ is a Virasoro minimal 
model and $S$ consists of one field $\phi_{2,1}(z)$, 
these questions have been studied in \cite{FJMMT1,FJMMT2}
under certain conditions. 
In this paper we consider the $\phi_{1,3}(z)$ field.
The corresponding filtration on the Virasoro modules
is called the $(1,3)$-filtration.
We also clarify the connection between 
the $\phi_{1,3}$ case and the 
fusion filtration on the representations of $\slth$.
Let us describe our results.

Let $p< p'$ be relatively prime positive integers, 
and let 
$M_{r,s}^{(p,p')}$ ($1\le r< p$, $1\le s<  p'$) 
be the irreducible representations of 
the Virasoro algebra with the central charge
$c=13-6(t+1/t)$ and highest weight 
$\triangle_{r,s}=((rt-s)^2-(t-1)^2)/4t$, 
where $t=p'/p$.  
We consider the $(1,3)$ primary field $(s-s'=-2,0,2)$
$$
\phi^{(s',s)}(z)=\sum_{n\in\Z+\triangle_{r,s'}- \triangle_{r,s}}
\phi_{-n}^{(s',s)} z^{n-\triangle_{1,3}},\
\quad
\phi_{-n}^{(s',s)}: M_{r,s}^{(p,p')}\to M_{r,s'}^{(p,p')}.
$$
We denote by $E_m(M_{r,s}^{(p,p')})$ the (1,3)-filtration
of $M_{r,s}^{(p,p')}$.
One of our goals is to construct 
a monomial basis of $M_{r,s}^{(p,p')}$compatible with
the (1,3) filtration under the condition $1< t< 2$. 
Fix a highest weight vector $|r,s\rangle$ of $M_{r,s}^{(p,p')}$.  
The proposed basis consists of vectors of the form 
\begin{equation}
\label{basis}
\phi_{-n_1}^{(s_0,s_1)}\ldots \phi_{-n_m}^{(s_{m-1},s_m)}|r,b(r,s)\rangle
\quad(s_0=s,s_m=b(r,s))
\end{equation}
where $b(r,s)$ is defined by
$$
\triangle_{r,b(r,s)}=\min_{1\le a\le p'-1\atop a\equiv s\bmod 2}
\triangle_{r,a}.
$$
The condition on the sequence $(n_1,\ldots, n_m)$ is
\begin{equation}
\label{cond}
n_i-n_{i+1}\ge w(s_{i-1},s_i,s_{i+1}),
\end{equation}
where $w(a,b,c)$ are some rational numbers
to be introduced in the text (see Section \ref{MONO}).  
We call vectors (\ref{basis}) subject to 
the condition (\ref{cond})
\emph{admissible monomials}.
Our conjecture is that the set 
of admissible monomials form a basis 
of $M_{r,s}^{(p,p')}$, 
which is compatible with the $(1,3)$-filtration. 
To support the conjecture, we prove two statements:
\begin{itemize}
\item[A.]
The character of the proposed basis
coincides with the character of $M_{r,s}^{(p,p')}$.
\item[B.]
In the unitary case $p'=p+1$, 
the $(q,v)$-character of the proposed basis 
and that of the $(1,3)$-filtration coincide.  
Here the $v$-degree of a monomial (\ref{basis}) 
is defined as the number $m$ of the operators $\phi_{-n_i}^{(s_{i-1},s_i)}$.
\end{itemize}

The proof of the statement A 
is based on combinatorics 
and the Rocha-Caridi character formula. 
Namely we show that the $q$-character of 
admissible monomials with fixed $m$ can be 
written in the form 
\begin{equation}
\label{Im}
\frac{q^{\triangle_{r,s}}}{(q)_m} I_m(q), 
\qquad 
(q)_m=\prod_{i=1}^m (1-q^i),
\end{equation}
where $I_m(q)$ is an alternating sum of
the characters for the fusion products \cite{FL}.
We also prove that the 
Rocha-Caridi formula for the character of 
$M_{r,s}^{(p,p')}$ can be rewritten as 
$\sum_{m\ge 0} \frac{q^{\triangle_{r,s}}}{(q)_m} I_m(q)$. 
In order to prove that 
the admissible monomials form a basis
of $M_{r,s}^{(p,p')}$,  
it is enough to rewrite any monomial 
$\phi_{-n_1}^{(s_0,s_1)}\cdots \phi_{-n_m}^{(s_{m-1},s_m)}|r,b(r,s)\rangle$
in terms of admissible monomials of 
length less than or equal to $m$. 
This was done in \cite{FJMMT1,FJMMT2} in the case of 
the $(2,1)$ field using 
quadratic relations for its Fourier components. 
Similar quadratic relations can be written also 
for the $(1,3)$ field using the results of \cite{DF}.  
Still it is not clear to us how to 
rewrite an arbitrary monomial 
in terms of admissible ones using these relations.

The proof of the statement B 
is based on the coset construction. 
For $i=0,1$ and integers $r,k$ with $0\le r-1\le k$, let 
\begin{equation}
\label{tenspr}
L_{i,1}\T L_{r-1,k}=
\bigoplus_{1\le s\le k+2\atop s:\ s+r+i \text{ even}}
M_{r,s}^{(k+2,k+3)}\T L_{s-1,k+1}
\end{equation}
be the decomposition of the tensor product of the irreducible highest weight
representations of $\slth$. Here $L_{j,k}$ denotes the level $k$ module
with the highest weight $j$ with respect to $\slt\T 1\hk\slth$.
We denote by $v_{j,k}$ the highest weight vector.
Using a result of \cite{L},  
we establish the connection  between the $(1,3)$
filtration and the fusion filtration on the left hand side of
(\ref{tenspr}). 
Namely, consider the action of the algebra $\slth$ 
on $L_{i,1}$ and the corresponding
$S$-filtration $G_m(L_{i,1})$, 
where $S=\{e(z),h(z),f(z)\}$,  
$x(z)=\sum (x\T t^n) z^{-n-1}$, 
and $\{e,h,f\}$ is the standard basis of $\slt$.  
We call this filtration the Poincar{\'e}-Birkhoff-Witt (PBW) filtration.
We show that
\begin{equation}
\label{PBWf}
U(\slth)\cdot (G_m(L_{i,1})\T v_{r-1,k})=
\bigoplus_{1\le s\le k+2\atop s:\ s+r+i \text{ even}}
E_m(M_{r,s}^{(k+2,k+3)})\T L_{s-1,k+1},
\end{equation}
where $E_m(M_{r,s}^{(k+2,k+3)})$ is the $(1,3)$ filtration. 
Thus the study of this filtration can 
be reduced to the study of the left hand side of 
(\ref{PBWf}). 

We recall that, for two cyclic $\mathfrak{g}$-modules $V_1$ and $V_2$ with
cyclic vectors $v_1$ and $v_2$, the fusion filtration on the tensor product 
$V_1(z_1)\T V_2(z_2)$ of evaluation
representations of $\mathfrak{g}\T\C[u]$ is defined by
$$
F_m(V_1(z_1)\T V_2(z_2))=
\mathrm{span}
\{(g^{(1)}\T u^{i_1})\cdots (g^{(s)}\T u^{i_s}) \cdot (v_1\T v_2)
\},
$$
where $g^{(i)}\in\mathfrak{g}$ and $i_1+\ldots + i_s\le m$.
One can easily show that
\be
F_m(L_{i,1}\T L_{l,k})=U(\slth)\cdot (G_m(L_{i,1})\T v_{l,k}).
\en

Using (\ref{PBWf}) we express the $(q,v)$-character of the
$(1,3)$-filtration via the $(q,v)$-character of the 
PBW-filtration on $L_{i,1}$.
We thus get the $(q,v)$-version (\ref{Im}) of the 
Rocha-Caridi formula as an
alternating sum of the $(q,v)$-characters of weight 
subspaces of $L_{i,1}$.

We finish the introduction with a discussion of possible
generalizations. Note that the integer level $k$ 
in the coset construction
(\ref{tenspr}) can be replaced by a fractional one. 
This generalization
leads to the coset realization of 
the general $M^{(p,p')}_{r,s}$. 
We expect
that the above construction can be applied to the general case.

Now consider the $(1,2)$ field $\phi_{1,2}(z)$. In this case the decomposition 
(\ref{tenspr}) should be replaced by
$$
(L_{0,1}\oplus L_{1,1})\T L_{r-1,k}=
\bigoplus_{1\le s\le k+2}
M_{r,s}^{(k+2,k+3)}\T L_{s-1,k+1}, 
$$
and the algebra $\slth$ by the vertex (intertwining) operator
$\C^2(z)$ acting on $L_{0,1}\oplus L_{1,1}$. This vertex operator induces
the filtration $G'_m(L_{0,1}\oplus L_{1,1})$. 
Then
$$
U(\slth)\cdot\left[G'_m(L_{0,1}\oplus L_{1,1})\T L_{r-1,k}\right]
=\bigoplus_{1\le s\le k+2} 
H'_m(M_{r,s}^{(k+2,k+3)})\T L_{s-1,k+1},
$$
where $H'_m(M_{r,s}^{(k+2,k+3)})$ is the $(1,2)$ filtration. 
As in the $(1,3)$-case, the $(q,v)$-character
$$
\sum_{m\ge 0} v^m \ch_qH'_m(M_{r,s}^{(k+2,k+3)})
$$ 
can be expressed as
an alternating sum of the $(q,v)$-characters of 
$\C^2(z)$-filtration on $L_{0,1}\oplus L_{1,1}$.

Our paper is organized as follows:
In Section 2, the character of the Virasoro module
$M^{(p,p')}_{r,s}$ is written as an alternating sum
by using the character of the weight $l$ component of the fusion product
$\pi_2^{*m}$ (Lemma \ref{lem:rocha}).
The main result in this section is a proof of the statement that
for $1< p'/p< 2$ the alternating sum with fixed $m$ is the character of
the admissible monomials of length $m$ (Proposition \ref{PROP1}).
Section 3 prepares some exact sequences of the fusion products and
vanishing of the homology groups of the Lie subalgebra
$\mathfrak n_+\subset\widehat{\mathfrak{sl}}_2$ generated by
$f\otimes1$ and $e\otimes t^{-1}$
with coefficients in the tensor product
of certain finite dimensional modules and irreducible highest weight
$\widehat{\mathfrak{sl}}_2$-modules. This is used in Section 5 when
the characters of the highest weight vectors in integrable 
$\widehat{\mathfrak{sl}}_2$-modules are computed.
In Section 4, by using Lashkevich's construction of vertex operators
in the GKO construction, an isomorphism is given between the fusion product
of level $1$ and level $k$ irreducible highest weight 
$\widehat{\mathfrak{sl}}_2$-modules and the associated graded module
with respect to the filtration defined by the $(1,3)$ primary field
(Proposition \ref{fusphi}). Section 5 is devoted to the calculation of
the characters for the $m$-th graded components of the Virasoro unitary series
with the $(1,3)$ filtration (Theorem \ref{MAIN}). The result coincides with the
combinatorial characters computed in Section 2.

Throughout the text, $e,f,h$ denotes the standard basis of $\slt$, and 
$\pi_j$ denotes the $(j+1)$-dimensional irreducible representation.
We use three different degree countings:
$q$-degree, $z$-degree and $v$-degree. The $q$-degree counts the eigenvalue
of the scaling operator $d$ in $\widehat{\mathfrak{sl}}_2$, or the power
in $u$ for $\mathfrak{sl}_2\otimes\C[u]$. The $z$-degree counts the eigenvalue
of the operator $h$ or $h\otimes1$. The $v$-degree counts the grading
by vertex operators $\phi_{1,3}(z)$ or $x(z)$ $(x\in\mathfrak{sl}_2)$.
\section{Conjectural monomial basis by (1,3) field}\label{MONO}

\subsection{Formulation}\label{subsec:formulation}

In this section, we consider 
Virasoro modules in the minimal series which are not necessarily unitary. 

Let ${\rm Vir}$ be the Virasoro algebra:
\be
[L_m,L_n]=(m-n)L_{m+n}+\frac c{12}m(m^2-1)\delta_{m+n,0}.
\en
Throughout this section, we fix relatively prime positive integers $p,p'$ 
satisfying $3\leq p< p'$. 
We denote by $M^{(p,p')}_{r,s}$ ($1\leq r\leq p-1$, $1\leq s\leq p'-1$) 
the irreducible {\rm Vir}-module with central charge $c=13-6(t+\frac1t)$  
and highest weight
\be
\Delta_{r,s}
=\frac{(rt-s)^2-(t-1)^2}{4t}\,, 
\en
where
\be
t=\frac{p'}{p}\,.
\en

We fix $1\leq r\leq p-1$, and consider the $(1,3)$ primary field
\be
\phi^{(s',s)}(z)=\sum_{n\in\Z+\Delta_{r,s'}-\Delta_{r,s}}
\phi^{(s',s)}_{-n}z^{n-\Delta_{1,3}}.
\en
The Fourier coefficients $\phi^{(s',s)}_{-n}$ are operators acting as
$\l(M^{(p,p')}_{r,s}\r)_d\rightarrow\l(M^{(p,p')}_{r,s'}\r)_{d+n}$, 
where 
$\Bigl(M^{(p,p')}_{r,s}\Bigr)_d=
\{|v\rangle\in M^{(p,p')}_{r,s}\mid L_0|v\rangle=d|v\rangle\}$
stands for the graded component. 
They are characterized by the intertwining property
\be
[L_n,\phi^{(s',s)}(z)]=z^n\l(z\frac d{dz}+(n+1)\Delta_{1,3}\r)\phi^{(s',s)}(z).
\en
A non-trivial $(1,3)$ primary field exists if and only if $s'=s,s\pm2$
and $(s',s)\not=(1,1),(p'-1,p'-1)$. Moreover, it is unique up to a
constant multiple. 
We fix the highest weight vector
$|r,s\rangle\in\Bigl(M^{(p,p')}_{r,s}\Bigr)_{\Delta_{r,s}}$
and use the normalization
$\phi^{(s',s)}_{\Delta_{r,s}-\Delta_{r,s'}}|r,s\rangle=|r,s'\rangle$.

Our problem is to construct a basis of the representation 
$M^{(p,p')}_{r,s}$ by using the operators $\phi^{(s',s)}_{-n}$.
In this paper we restrict to the case 
\bea
1< t< 2, \label{eq:1< t< 2}
\ena
and give a partial answer to this problem. 

The form of the basis we propose is similar to the one studied in 
\cite{FJMMT1,FJMMT2} using the $(2,1)$ primary field. 
We define a set of weights
\bel
\label{WEI}
w(a,b,c) \in \Delta_{r,a}-2\Delta_{r,b}+\Delta_{r,c}+{\Z}\,,
\enl
and consider vectors of the form
\bel
\label{MON}
\phi^{(s_0,s_1)}_{-n_1}\cdots\phi^{(s_{m-1},s_m)}_{-n_m}|r,s_m\rangle
\enl
satisfying 
\bel
\label{ADM}
n_i-n_{i+1}\geq w(s_{i-1},s_i,s_{i+1})\,.
\enl
The actual form of the weights \eqref{WEI} is a little involved, 
and we postpone their definition to subsection \ref{subsec:weight} 
(see \eqref{W1}--\eqref{W6} and paragraphs following them).   
The vectors \eqref{MON} 
are parametrized by spin-$1$ and level-$p'$ restricted paths,
i.e., sequences of integers
\be
{\bf s}=(s_0,s_1,\ldots,s_m)
\en
satisfying
\be
& & 1\leq s_i\leq p'-1,\\
& & s_i=s_{i+1}\hbox{ or }s_{i+1}\pm2,\\
& & (s_i,s_{i+1})\not=(1,1),(p'-1,p'-1).
\en
We call them simply {\it paths}. 
The non-negative integer $m$ is called the length of the path. 
We denote by $\P^{(p')}_{a,b,m}$ the set of paths 
$(s_0,s_1,\ldots,s_m)$ of length $m$ satisfying $s_0=a$ and $s_m=b$.
Note that the parity of $s_i$ is common with each path. 
In particular, we have $a\equiv b\bmod2$ if $\P^{(p')}_{a,b,m}$ is non-empty.

The set of rational numbers
\be 
(n_1,\ldots,n_m)
\en
in the expression (\ref{MON}) is called a rigging associated with the path 
$(s_0,s_1,\ldots,s_m)$. A rigging satisfies
\be
n_i\in\Z+\Delta_{r,s_{i-1}}-\Delta_{r,s_i}.
\en
A path with rigging is called a rigged path.

A rigged path of length $m$ is called 
\emph{admissible} 
if and only if
(\ref{ADM}) and the following condition hold:
\bel\label{BDR}
n_m\geq\Delta_{r,s_{m-1}}-\Delta_{r,s_m}+\delta_{s_{m-1},s_m}.
\enl
We denote by $\R^{(p,p')}_{r,a,b,m}$ the set of admissible rigged paths
of length $m$ such that $s_0=s$ and $s_m=b$. 

Finally, we fix the boundary of a path $s_m=b(r,s)$ 
by the following rule:  $b=b(r,s)$ is the unique integer
satisfying $1\le b\le p'-1$, $b\equiv a\bmod 2$ and 
\bea
\Delta_{r,b}=
\min_{1\le s\le p'-1\atop s\equiv a\bmod 2}
\Delta_{r,s}\,.
\label{RAB}
\ena
We will comment on this choice in the next subsection \ref{subsec:rewriting}.

Now we put forward the 

\begin{conj}\label{CONJ}
The set of vectors \eqref{MON},  
where $(n_1,\ldots,n_m)$ runs through
the set $\bigcup_{m\ge 0}\R^{(p,p')}_{r,a,b,m}$ and 
$b=b(r,a)$ as given in \eqref{RAB}, 
constitute a basis of $M^{(p,p')}_{r,a}$.
\end{conj}

Note that the meaning of the condition (\ref{BDR}) is clear.
If $\Delta_{r,s_m}+n_m< \Delta_{r,s_{m-1}}$, the vector
$\phi^{(s_{m-1},s_m)}_{-n_m}|r,s_m\rangle$ is zero because
$\Bigl(M^{(p,p')}_{r,s_{m-1}}\Bigr)_{\Delta_{r,s_m}+n_m}=\{0\}$.
If $s_{m-1}=s_m$ and $n_m=0$, the vector 
$\phi^{(s_m,s_m)}_0|r,s_m\rangle$ is proportional to
$|r,s_m\rangle$ because it belongs to 
$\Bigl(M^{(p,p')}_{r,s_m}\Bigr)_{\Delta_{r,s_m}}=\C|r,s_m\rangle$.
We have also $\Delta_{1,1}=0$ and
$\Bigl(M^{(p,p')}_{1,1}\Bigr)_1=\{0\}$. 
This is not in contradiction to
the condition (\ref{BDR}) for $r=s_m=s_{m-1}=1$
because the case $s_m=s_{m-1}=1$ is prohibited.

In order to support the conjecture, we will show that 
the set of admissible monomials has the same character 
as that of $M^{(p,p')}_{r,a}$ (Theorem \ref{PROP2}). 
The conjecture will follow if we show further that the above 
set of vectors span the space $M^{(p,p')}_{r,a}$.
So far we have not been able to check the latter point in full generality. 

\subsection{A character identity}\label{subsec:rewriting}
In this subsection, we rewrite the character of $M^{(p,p')}_{r,a}$ in a form 
suitable for comparison with the set of paths. 

We use the $q$-supernomial coefficients introduced in \cite{SW}.  
They are a $q$-analog of the weight multiplicities of tensor products of 
various $\pi_k$, where $\pi_k$ denotes the 
irreducible $\slt$-module of dimension $k+1$. 
As shown in \cite{FF1}, they can be  
defined as the coefficients of $z^l$ of 
the graded character of the fusion product 
(for the definition and properties of fusion product, see Section 3).
\be
{\rm ch}_{q,z}\Bigl(
\pi_1^{*L_1}*\cdots*\pi_N^{*L_N}
\Bigr)
=
\sum_{l\in\Z+\frac{1}{2}\sum_{j=1}^NjL_j}
\qbin{L_1,\cdots,L_N}{l\phantom{\cdots,L_N}}_qz^l \,.
\en
Here we will need only the special case $N=2$, $L_1=0$.  
Set 
\bea
S_{m,l}(q):=\qbin{0,m}{l\phantom{ L_2 }}_{q^{-1}}.
\label{eq:Sml}
\ena
Formula (2.9) in \cite{SW} gives
\be
S_{m,l}(q)=\sum_{\nu\in\Z}
q^{(\nu+l-m)(\nu+l)+\nu(\nu-m)}
{\qbin{m}{\nu}}_q{\qbin{\nu}{m-l-\nu}}_q.
\en
In the right hand side 
\be
{\qbin{L}{a}}_q=
\begin{cases}
\frac{(q^{L-a+1})_a}{(q)_a} &  (a\in \Z_{\ge 0},L\in \Z),\\
0 & (\mbox{otherwise}),\\
\end{cases}
\en
stands for the $q$-binomial symbol,  
and $(x)_n=(x)_\infty/(q^nx)_\infty$, 
$(x)_\infty=\prod_{i=0}^\infty(1-q^ix)$. 

Recall that the character $\chi^{(p,p')}_{r,s}(q)$ of $M^{(p,p')}_{r,s}$ 
is given by \cite{RC}
\bea
q^{-\Delta^{(p,p')}_{r,s}}
\chi^{(p,p')}_{r,s}(q)
& = &
\sum_{\la\in\Z}
\frac{1}{(q)_\infty}
q^{\la^2 p p'+\la(p' r-ps)}
- -\sum_{\la\in\Z}
\frac{1}{(q)_\infty}
q^{\la^2 p p'+\la(p' r+ps)+rs}
\,.
\label{eq:rocha}
\ena

\begin{lem}\label{lem:rocha}
For any $b\in\Z$ satisfying $b\equiv a\bmod 2$,  
the character \eqref{eq:rocha} 
can be written in terms of \eqref{eq:Sml} as 
\bea
& & q^{-\Delta^{(p,p')}_{r,a}}\chi^{(p,p')}_{r,a}(q)
=\sum_{m\ge 0}\frac{1}{(q)_m}I^{(p,p')}_{r,a,b,m}(q),
\label{eq:rocha2}
\ena
where
\bea
I^{(p,p')}_{r,a,b,m}(q)
& =& \sum_{\la\in\Z}q^{\la^2 p p'+\la(p' r-pa)+m^2-((a-b)/2-p'\la)^2}
S_{m,(a-b)/2-p'\la}(q)
\label{eq:Im}
\\
& - &
\sum_{\la\in\Z}q^{\la^2 p p'+\la(p' r+pa)+ra+m^2-((a+b)/2+p'\la)^2}
S_{m,(a+b)/2+p'\la}(q)
\,.
\nn
\ena
\end{lem}
\begin{proof}
The $q$-supernomial coefficients satisfy the recurrence relations
(\cite{SW}, Lemma 2.3) 
\be
& & 
\qbin{L_1,L_2}{a\phantom{ L_2 }}_q
=q^{L_1+L_2-1}\qbin{L_1-2,L_2}{a\phantom{ L_2 }}_q
+\qbin{L_1-2,L_2+1}{a\phantom{ L_2 }}_q\,
\\
& & \qbin{L_1,0}{a\phantom{ L_2 }}_q=
{\qbin{L_1}{a+L_1/2}}_q.
\en
Iterating this $k$ times, we find 
\be
& & \qbin{L_1,L_2}{a\phantom{ L_2 }}_q
=\sum_{m=0}^kq^{(k-m)(L_1+L_2-k)}{\qbin{k}{m}}_q
\qbin{L_1-2k,L_2+m}{a\phantom{ L_2 }}_q\,.
\en
Choosing $L_1=2N,L_2=0,k=N$, changing $q\to q^{-1}$ and 
letting $N\to\infty$ we obtain for all $l\in \Z$ 
\be
\frac{1}{(q)_\infty}=\sum_{m\ge 0}\frac{q^{m^2-l^2}}{(q)_m}
S_{m,l}(q)\,.
\en

In each summand of the the first (resp. second) sum
of \eqref{eq:rocha}, replace 
$1/(q)_\infty$ by the right hand side
of the above identity, 
choosing $l=(a-b)/2-p'\la$ (resp. $(a+b)/2+p'\la$). 
The desired identity follows.
\end{proof}

Though \eqref{eq:rocha2} is an identity valid for any $b\in\Z$, 
in most cases the polynomial \eqref{eq:Im} comprises negative coefficients. 
We prove in subsection \ref{subsec:proof} that,  
if $1< t< 2$ and $(r,b)$ satisfies (\ref{RAB}), then 
the coefficients of $I^{(p,p')}_{r,a,b,m}(q)$ are non-negative integers.
In fact we will show that it can be written 
as a configuration sum over the set of paths $\P^{(p')}_{a,b,m}$. 
Define the weight of a path ${\bf s}\in\P^{(p')}_{a,b,m}$ by
\be
E({\bf s})=\sum_{i=1}^{m-1}iw(s_{i-1},s_i,s_{i+1}).
\en

\begin{prop}\label{PROP1}
Under the conditions $1< t< 2$ and $(\ref{RAB})$, 
we have an equality
\bel\label{GEN}
I^{(p,p')}_{r,a,b,m}(q)=\sum_{{\bf s}\in\P^{(p')}_{a,b,m}}
q^{E({\bf s})+m(\Delta_{r,s_{m-1}}-\Delta_{r.b}+\delta_{s_{m-1},b})
+\Delta_{r,b}-\Delta_{r,a}}\,.
\enl
\end{prop}
We give a proof in Section 2.4.
Note that the exponent of $q$ in (\ref{GEN})
is an integer because of (\ref{WEI}).

{}From Proposition \ref{PROP1} immediately follows
\begin{thm}\label{PROP2} Notation being as above, 
we have an identity for the character 
\be
\chi^{(p,p')}_{r,a}(q)=
\sum_{m\geq0}
\frac1{(q)_m}
\sum_{{\bf s}\in\P^{(p')}_{a,b,m}}
q^{E({\bf s})+m(\Delta_{r,s_{m-1}}-\Delta_{r.b}+\delta_{s_{m-1},b})
+\Delta_{r,b}}.
\en
\end{thm}

\subsection{Definition of $w(a,b,c)$}\label{subsec:weight}

In this subsection we introduce our weight $w(a,b,c)$. 

In \eqref{GEN}, we fixed $b=b(r,a)$ by the condition (\ref{RAB}).
Conversely, for a given $b$, $r$ for which (\ref{RAB}) is valid
is either $r_1(b)=\left[\frac{b+1}t\right]$ or 
$r_2(b)=\left[\frac{b-1}t\right]+1$.
Set
\be
\tau(b)=\left[\frac{b+1}t\right]-\left[\frac{b-1}t\right].
\en
We have $\tau(b)=1$ or $2$. If $\tau(b)=1$ we have $r_1(b)=r_2(b)$,
and if $\tau(b)=2$ we have $r_1(b)=r_2(b)+1$.
We list a few other properties of $\tau(s)$.
\bel
& & \tau(1)=1,\label{T1}\\
& & \tau(2)=\begin{cases}
2& \hbox{if }1< t< \frac32;\\
1& \hbox{if }\frac32< t< 2,
\end{cases}\label{T2}\\
& & \tau(p'-1)=2,\label{T3}\\
& & \tau(s)=\tau(p'-s)\hbox{ if }1< s< p'-1.\label{T4}
\enl

Set
\be
\{x\}=x-[x]
\en
where $[x]$ is the integer part of $x$.

We define the weight $w(a,b,c)$ in the following form:
\bel
& & w(s\pm2,s,s\mp2)=\frac2t,\label{W1}\\
& & w(s,s+2,s+2)=w(s+2,s+2,s)=2-\left\{\frac{s+1}t\right\},\label{W2}\\
& & w(s,s,s+2)=w(s+2,s,s)=1+\left\{\frac{s+1}t\right\},\label{W3}\\
& & w(s,s+2,s)=-2\left\{\frac{s+1}t\right\}+x(s),\label{W4}\\
& & w(s,s,s)=3-\tau(s),\label{W5}\\
& & w(s,s-2,s)=2\left\{\frac{s+1}t\right\}-\frac4t+y(s).\label{W6}
\enl
Here $x(s),y(s)$ are integers given by one of the following:
\be
\begin{matrix}
               &  \tau(s) &  x(s) &  y(s) \\
\mbox{Case $1_A$} &   1      &  2    &   3   \\
\mbox{Case $1_B$} &   1      &  3    &   2   \\
\mbox{Case 2}  &  2      &  2    &   4   \\
\end{matrix}
\en
If $\tau(s)=2$, then the choice of $x(s),y(s)$ is unique. 
Otherwise there are two cases. We write 
$\tau(s)=1_A$ or $\tau(s)=1_B$ to indicate which of the two rules apply. 

If $1< t\le 5/3$ and $\tau(s)=1$, then we take
\be
& & 
\tau(s)=\begin{cases}
1_A &  (1\le s< \frac{p'}{2}),\\
1_B &  (\frac{p'}{2}< s \le p'-1).\\
\end{cases}
\en
If $5/3< t< 2$, then we take 
\be
& & \tau(1),\tau(2),\cdots,\tau(p'-1)
\\
& & =1_A,1_B,1_A,1_B,\ldots,1_B,1_A,2,1_B,1_A,1_B,\ldots,1_B,1_A,2,
\\
& & \quad 1_B,1_A,1_B,\ldots, 1_A,2, 1_B,1_A,1_B,\ldots,1_B,1_A,2,
\en
where $1$'s between successive $2$'s come always with an even number. 
Below we will motivate the above assignment and show 
that it can be made consistently.

Let us seek for the weights $w(a,b,c)$ in the above form (\ref{W1}-\ref{W6}). 
We will take them independently of the choice $r=r_1(b)$ or $r_2(b)$. 
We demand further the following. 
\begin{enumerate}
\item (\ref{GEN}) holds for $m=2$, 
\item left-right symmetry $w(a,b,c)=w(c,b,a)$,  
\item symmetry reflecting $M^{(p,p')}_{r,s}\simeq M^{(p,p')}_{p-r,p'-s}$, 
\bel\label{SYM}
w(a,b,c)=w(p'-a,p'-b,p'-c)\,,
\enl
\item 
\bel\label{R3}
\hbox{if }\tau(s)=1_B\hbox{ then }\tau(s+2)\not=1_A.
\enl
\end{enumerate}
The last condition turns out to be necessary in the course of the proof
of (\ref{GEN}), see subsection \ref{subsec:proof}. 

The validity of \eqref{GEN} for $m=2$ gives a linear constraint on
the weights. 
There are three cases: $a=b\pm4,b\pm2,b$.
In the first two cases, the relevant weights are (\ref{W1}-\ref{W3}).
They are independent of $\tau(s)$. It is easy to check that
the constraint is satisfied in these cases.
In the third case, the relevant weights are (\ref{W4}-\ref{W6}).
Here the value of $\tau(s)$ matters. If $\tau(s)=2$, the weights
are uniquely given by (\ref{W4}-\ref{W6}), and they satisfy the constraint.
If $\tau(s)=1$, we must specify $\tau(s)=1_A$ or $1_B$.
For $s=1,2$, see (\ref{T1}) and (\ref{T2}). In these cases,
the constraint implies 
\bel
\tau(1)& =& 1_A,\label{C1}\\
\tau(2)& =& 1_B\quad\hbox{for}\quad\frac32< t< 2.\label{C2}
\enl
If $2< s< p'-2$, the constraint is satisfied for either choice.

The left-right symmetry (ii) is automatically satisfied by 
the formulas (\ref{W1}-\ref{W6}). 

Symmetry \eqref{SYM} is also valid for (\ref{W1}-\ref{W3}) and (\ref{W5}). 
it is obvious for (\ref{W1}); and follows from
\be
\left\{\frac st\right\}+\left\{-\frac st\right\}=1\hbox{ if } 1< s< p'
\en
for (\ref{W2}) $\leftrightarrow$ (\ref{W3}), and for (\ref{W5}). 

We determine the choice of $1_{A/B}$ so that 
the symmetry (\ref{SYM}) for (\ref{W4}) $\leftrightarrow$ (\ref{W6}) is valid.
Since $\tau(1)=1_A$ and $\tau(p'-1)=2$, we have
\be
w(1,3,1)=w(p'-1,p'-2,p'-1)=4-\frac4t.
\en
For $1< t< 3/2$, we have $\tau(2)=\tau(p'-2)=2$, and
\be
w(2,4,2)=w(p'-2,p'-4,p'-2)=6-\frac6t;
\en
for $3/2< t< 2$, we have $\tau(2)=1_B$. Setting $\tau(p'-2)=1_A$, we have
\be
w(2,4,2)=w(p'-2,p'-4,p'-2)=5-\frac6t.
\en
For $2< s< p'-2$, the symmetry is valid if $\tau(s)=\tau(p'-s)=2$;
if $\tau(s)=\tau(p'-s)=1$, we need to choose $1_{A/B}$ in such a way that
\bel\label{R2}
\tau(s)=1_A,1_B\leftrightarrow\tau(p'-s)=1_B,1_A.
\enl
The last requirement would be inconsistent if $\tau(s)=1$ for $s=p'/2$.
However, we have
\begin{lem}
If $p'$ is even, we have $\tau(p'/2)=2$.
\end{lem}
\begin{proof}
If $p'$ is even, then $p$ is odd.
We have
\be
\tau(p'/2)=[\alpha]-[\beta]\hbox{ where }
\alpha=\frac p2+\frac p{p'},\ \beta=\frac p2-\frac p{p'}.
\en
Since $\alpha+\beta=p$ is odd and $1< \alpha-\beta=2p/p'< 2$,
we have $\tau(p'/2)=2$.
\end{proof}

We need also to satisfy the condition \eqref{R3}. 
Let us show the consistency of (\ref{C1}), (\ref{C2}), (\ref{R2}) and (\ref{R3}).
Suppose that $\tau(s)=1$. We choose $\tau(s)=1_A$ or $1_B$ as follows.

If $1< t< 3/2$, we have $\tau(1)=1_A$ and $\tau(2)=2$.
Therefore, the following choice of $1_A$ or $1_B$ for $s$ such that $\tau(s)=1$
satisfies all the constraints.
\bel\label{AB1}
\tau(s)=
\begin{cases}
1_A& \hbox{ if }s< p'/2;\\
1_B& \hbox{ if }s> p'/2.
\end{cases}
\enl

If $3/2< t\leq5/3$, we have $\tau(1)=1_A$, $\tau(2)=1_B$ and $\tau(4)=2$.
Therefore, the same choice (\ref{AB1}) will do.

Before going to the case $5/3< t< 2$, we prepare a few lemmas.

\begin{lem}\label{L1}
If $3/2< t< 2$, we do not have the sequence $(\tau(s),\tau(s+1))=(2,2)$.
\end{lem}
\begin{proof}
Since $\frac12< \frac1t< 1$, the increment $[\frac{s+1}t]-[\frac st]$
is either $0$ or $1$. Therefore, if $\tau(s)=\tau(s+1)=2$,
we have
\be
& & m\leq \frac{s-1}t< m+1,\quad m+1\leq \frac st< m+2,\\
& & m+2\leq \frac{s+1}t< m+3,\quad m+3\leq \frac{s+2}t< m+4
\en
for some integer $m$. From $\frac{s-1}t< m+1$ and $m+3\leq \frac{s+2}t$
follows $t< \frac32$.
\end{proof}
\begin{lem}\label{L2}
Suppose that $\left[\frac st\right]=m$,
$\left[\frac{s+1}t\right]=m$. Then, $\left[\frac{s+2}t\right]=m+1$.
\end{lem}
\begin{proof}
The statement follows from $1< \frac2t< 2$.
\end{proof}
\begin{lem}\label{L3}
Suppose that $\tau(s)=2$, $\tau(s+1)=\cdots=\tau(s+k)=1$ and
$\tau(s+k+1)=2$. Then, $k$ is even.
\end{lem}
\begin{proof}
By Lemma \ref{L2}, we have the sequence
\be
& & {\textstyle\left[\frac{s-1}t\right]}=m,
{\textstyle\left[\frac st\right]}=m+1,
{\textstyle\left[\frac{s+1}t\right]}=m+2,
{\textstyle\left[\frac{s+2}t\right]}=m+2,\\
& & {\textstyle\left[\frac{s+3}t\right]}=m+3,
{\textstyle\left[\frac{s+4}t\right]}=m+3,\ldots
{\textstyle\left[\frac{s+k-1}t\right]}=m+l,
{\textstyle\left[\frac{s+k}t\right]}=m+l,\\
& & {\textstyle\left[\frac{s+k+1}t\right]}=m+l+1,
{\textstyle\left[\frac{s+k+2}t\right]}=m+l+2.
\en
Therefore, $k=2l$.
\end{proof}

If $5/3< s< 2$, we have $\tau(1)=1_A$ and $\tau(2)=1_B$.
If we determine the choice for $s< p'/2$, the rest is determined
by (\ref{R2}). The constraint (\ref{R3}) together with the symmetry
(\ref{R2}) implies that if $s+2< p'/2$ and $\tau(s)=\tau(s+2)=1$ we have
\be
\tau(s)=1_A,1_B\rightarrow\tau(s+2)=1_A,1_B.
\en
We start from $\tau(1)=1_A$ and $\tau(2)=1_B$
and continue as $1_A,1_B,1_A,1_B,\ldots$ until $2$ appears.
By a similar argument as in the proof of Lemma \ref{L3},
we can show that the first appearance of $2$ is for even $s$.
Therefore, from Lemmas \ref{L1} and \ref{L3} we can define the sequence
$\tau(1),\tau(2),\ldots$ as
\be
1_A,1_B,1_A,1_B,\ldots,1_B,1_A,2,1_B,1_A,1_B,\ldots,
1_B,1_A,2,1_B,1_A,1_B,\ldots
\en
This sequence does not contain $(\tau(s),\tau(s+2))=(1_B,1_A)$.
The constraint (\ref{R3}) is also satisfied.

\subsection{Proof of Proposition \ref{PROP1}}\label{subsec:proof}

In this subsection we fix the weights $w(a,b,c)$ as in the previous
section, and prove \eqref{GEN}. 
To that end we consider the configuration sum
\begin{eqnarray*}
X_{a,b,c,m}(q):=\sum_{{\bf s} \in 
\mathcal{P}_{a,b,c,m}^{(p')}}q^{E({\bf s})}, 
\end{eqnarray*}
where $\mathcal{P}_{a,b,c,m}^{(p')}$ is 
the set of paths $(s_{0}, \ldots , s_{m+1})$
satisfying $s_{0}=a, s_{m}=b$ and $s_{m+1}=c$. 
{}From the definition, 
$X_{a,b,c,m}(q)=0$ unless $1\le a,b,c \le p'-1$, 
$b=c,c\pm 2$ and $(b,c)\neq (1,1),(p'-1,p'-1)$.
Note that $X_{a,b,c,m}(q)$ is uniquely determined from 
the initial condition $X_{a,b,c,0}(q)=\delta_{a,b}$ and 
the recurrence relation
\begin{eqnarray}
X_{a,b,c,m+1}(q)=\sum_{d=b, \, b\pm2}
q^{(m+1)w(d,b,c)}X_{a,d,b,m}(q). 
\label{eq:recurrenceX}
\end{eqnarray}

Let us give an explicit formula for $X_{a,b,c,m}(q)$. 
As an ingredient we introduce the function $\widetilde{S}_{m,l}(q)$ 
defined by 
\begin{eqnarray*}
\widetilde{S}_{m,l}(q):=\sum_{\nu \in \mathbb{Z}}
q^{(\nu+l-m)(\nu+l-1)+\nu(\nu-m)}
{\qbin{m}{\nu}}_q{\qbin{\nu}{m-l-\nu}}_q. 
\end{eqnarray*}
The functions $S_{m,l}(q)$ and $\widetilde{S}_{m,l}(q)$ 
are related to each other as follows. 

\begin{lem}
The following formulae hold:
\begin{eqnarray}
& &  
S_{m, -l}(q)=S_{m,l}(q), \qquad 
\widetilde{S}_{m, -l}(q)=q^{l}\,\widetilde{S}_{m,l}(q), 
\label{eq:relS0}\\ 
& &  
S_{m+1,l}(q)=
q^{-m-l-1}S_{m,l+1}(q)+S_{m,l}(q)+q^{-m+l-1}\widetilde{S}_{m,l-1}(q)\,,
\label{eq:relS1}\\ 
& &  S_{m+1,l}(q)=
q^{-m-l-1}S_{m,l+1}(q)+q^{-m}\widetilde{S}_{m,l}(q)+S_{m,l-1}(q)\,,
\label{eq:relS2}\\ 
& &  S_{m+1,l}(q)=
q^{-m}\widetilde{S}_{m,l+1}(q)+S_{m,l}(q)+q^{-m+l-1}S_{m,l-1}(q)\,, 
\label{eq:relS3}\\
& &  
\widetilde{S}_{m+1,l}(q)=
q^{-l}S_{m,l+1}(q)+S_{m,l}(q)+q^{-m+l-1}\widetilde{S}_{m,l-1}(q)\,, 
\label{eq:relS4}
\\ 
& &  \widetilde{S}_{m+1,l}(q)=
q^{-l}S_{m,l+1}(q)+q^{-m}\widetilde{S}_{m,l}(q)+S_{m,l-1}(q)\,,
\label{eq:relS5}\\ 
& &  \widetilde{S}_{m+1,l}(q)=
q^{-m-l}\widetilde{S}_{m,l+1}(q)+q^{-l}S_{m,l}(q)+S_{m,l-1}(q)\,.
\label{eq:relS6}
\end{eqnarray} 
\end{lem}
\begin{proof}
We use the notation of $q$-trinomial 
\begin{eqnarray*}
\qtri{n}{a}{b}{c}:=
\frac{(q)_{n}}{(q)_{a}(q)_{b}(q)_{c}}\qquad 
\hbox{for}\quad a+b+c=n. 
\end{eqnarray*}
Then the product of $q$-binomials in 
the definition of $S_{m,l}$ and $\widetilde{S}_{m,l}$ 
is rewritten as 
\begin{eqnarray*}
\qbin{m}{\nu}\qbin{\nu}{m-\nu-l}=
\qtri{m}{m-\nu}{m-l-\nu}{2\nu+l-m} 
\end{eqnarray*} 
{}From this expression it is easy to check \eqref{eq:relS0}. 
The other formulae except \eqref{eq:relS2} and 
\eqref{eq:relS5} follow directly from the $q$-trinomial identity: 
\begin{eqnarray}
\hspace{2em}
\qtri{n}{a}{b}{c}=\qtri{n-1}{a-1}{b}{c}+q^{a}\qtri{n-1}{a}{b-1}{c}+
q^{a+b}\qtri{n-1}{a}{b}{c-1}. 
\label{eq:qtri}
\end{eqnarray} 

In the following we prove \eqref{eq:relS2}. 
The proof of \eqref{eq:relS5} is similar. 

We start from 
\begin{eqnarray*}
S_{m+1,l}(q)=\sum_{\nu \in \mathbb{Z}}
q^{(\nu+l-m-1)(\nu+l)+\nu(\nu-m-1)}
\qtri{m+1}{m+1-l-\nu}{2\nu+l-m-1}{m+1-\nu}. 
\end{eqnarray*}
Decompose the right hand side above into three parts 
by applying \eqref{eq:qtri} to the $q$-trinomial.
Then by changing $\nu \to \nu+1$ we see that 
the first part is equal to $q^{-m-l-1}S_{m,l+1}(q)$. 
In the third part we rewrite the $q$-trinomial as follows: 
\begin{eqnarray*}
& &  
\qtri{m}{m+1-l-\nu}{2\nu+l-m-1}{m-\nu} \\ 
& &  {}=
\left(q^{m+1-l-\nu}+(1-q^{m+1-l-\nu})\right)
\qtri{m}{m+1-l-\nu}{2\nu+l-m-1}{m-\nu} \\ 
& &  {}=
q^{m+1-l-\nu}\qtri{m}{m+1-l-\nu}{2\nu+l-m-1}{m-\nu} \\ 
& &  \qquad {}+
(1-q^{2\nu+l-m})
\qtri{m}{m+l-\nu}{2\nu+l-m}{m-\nu}.
\end{eqnarray*}
{}From the first term in the right hand side above, 
we obtain $q^{-m}\widetilde{S}_{m,l}(q)$. 
Then the remaining is 
\begin{eqnarray*}
& &  
\sum_{\nu \in \mathbb{Z}}
q^{(\nu+l-m-1)(\nu+l)+\nu(\nu-m-1)+m+1-l-\nu}
\qtri{m}{m+1-l-\nu}{2\nu+l-m-2}{m-\nu+1} \\ 
& &  {}+
\sum_{\nu \in \mathbb{Z}}
q^{(\nu+l-m-1)(\nu+l)+\nu(\nu-m-1)+\nu}
(1-q^{2\nu+l-m})
\qtri{m}{m+l-\nu}{2\nu+l-m}{m-\nu}.
\end{eqnarray*}
Change $\nu \to \nu+1$ in the first sum, then 
it is canceled by the part 
containing $-q^{2\nu+l-m}$ in the second sum. 
The rest is equal to $S_{m, l-1}(q)$. 
\end{proof}

Now we define the function $f_{a,b,c,m}(q)$ 
in the following form 
for $(a,b,c)\in\mathbb{Z}^{3}$ satisfying 
$1 \le b, c \le p'-1, \, a \equiv b \, (\hbox{mod}\,2)$ and $c=b, \, b\pm2$:
\begin{eqnarray*}
& &  
f_{a,b,b+2,m}(q):=
q^{l(l+1)/t+m^2-l^2+(m-l)\{(b+1)/t\}}\times 
\left\{ 
\begin{array}{ll}
\widetilde{S}_{m,l}(q) &  \hbox{if}\quad \tau(b+2)=1_{A}, \\
S_{m,l}(q) &  \hbox{if}\quad \tau(b+2)=1_{B}\,\,\hbox{or}\,\,2, 
\end{array}
\right. \\
& &  
f_{a,b,b,m}(q):=
q^{l(l-1)/t+m^2-l^2+l(1-\{(b-1)/t\})}\times 
\left\{ 
\begin{array}{ll}
q^{m}S_{m,l}(q) &  \hbox{if}\quad \tau(b)=1_{A}\,\,\hbox{or}\,\,1_{B}, \\
q^{l}\widetilde{S}_{m,l}(q) &  \hbox{if}\quad \tau(b)=2, 
\end{array}
\right. \\
& &  
f_{a,b,b-2,m}(q):=
q^{l(l-1)/t+m^2-l^2+(m+l)(1-\{(b-1)/t\})}\times 
\left\{ 
\begin{array}{ll}
S_{m,l}(q) &  \hbox{if}\quad \tau(b-2)=1_{A}\,\,\hbox{or}\,\,2, \\
q^{l}\widetilde{S}_{m,l}(q) &  \hbox{if}\quad \tau(b-2)=1_{B}.
\end{array}
\right. 
\end{eqnarray*}
Here we set $l:=(b-a)/2$ in the right hand sides. 
By definition we set $f_{a,b,c,m}(q)=0$ for other type of 
triples $(a,b,c)$. 
Then $X_{a,b,c,m}$ is given in terms of $f_{a,b,c,m}$:

\begin{prop}\label{thm:Xandf}
We have 
\begin{eqnarray}
X_{a,b,c,m}(q)=\sum_{\epsilon=\pm}\epsilon
\sum_{n \in \mathbb{Z}}f_{\epsilon(a+2p'n), b, c, m}(q).
\label{eq:Xandf}
\end{eqnarray} 
\end{prop}
\begin{proof}
It is easy to check that the right hand side of \eqref{eq:Xandf} 
satisfies the initial condition, that is, 
it is equal to $\delta_{a,b}$ in the case of $m=0$. 
Now we should check the recurrence relation 
\begin{eqnarray}
\sum_{d=b, \, b\pm2}
q^{(m+1)w(d,b,c)}\sum_{\epsilon=\pm}\epsilon
\sum_{n \in \mathbb{Z}} 
f_{\epsilon(a+2p'n),d,b,m}(q)=
\sum_{\epsilon=\pm}\epsilon
\sum_{n \in \mathbb{Z}} 
f_{\epsilon(a+2p'n),b,c,m+1}(q). 
\label{eq:recurrencef}
\end{eqnarray}
Divide the cases according to the value of $b$\,;
(i) $3 \le b \le p'-3$ (non-boundary), 
(ii) $b=2$ or $p'-2$ (next-to-boundary) and 
(iii) $b=1$ or $p'-1$ (boundary). 
Then the proof of \eqref{eq:recurrencef} is 
just case-checking for each combination of the values 
$c\,(=b, b\pm2), \tau(b)$ and $\tau(c)$. 
In the following we give a sketch of the calculation 
in two cases as an example. 

First let us consider one of the non-boundary cases; 
$3 \le b \le p'-3, c=b-2$ and $\tau(b)=1_{A}$. 
Then from the definition of $f_{a,b,c,m}$ we have 
\begin{eqnarray*}
& &  
\sum_{d=b, b\pm2}q^{(m+1)w(d,b,b-2)}f_{a+2p'n, d, b, m}(q) \\ 
& &  \qquad {}=
q^{(l-p'n)(l-p'n-1)/t+(m+1)^{2}-(l-p'n)^{2}+(m+1+l-p'n)(1-\{(b-1)/t\})}
\\ 
& &  \qquad {}\times 
\Bigl\{ 
q^{-m-l+p'n-1}S_{m,l-p'n+1}(q)+S_{m,l-p'n}(q)+
q^{(m+1)(x(b-2)-3)+l-p'n}\widetilde{S}_{m, l-p'n-1}(q)
\Bigr\}. 
\end{eqnarray*} 
Here we set $l=(b-a)/2$.  
Since $\tau(b)=1_{A}$, the value of $\tau(b-2)$ is either 
$1_{A}$ or $2$. 
In both cases we have $x(b-2)=2$, and hence we can apply 
\eqref{eq:relS1} to the right hand side above. 
Thus we obtain 
\begin{eqnarray*}
\sum_{d=b, b\pm2}q^{(m+1)w(d,b,b-2)}f_{a+2p'n, d, b, m}(q)
 =f_{a+2p'n,b,b-2,m+1}(q). 
\end{eqnarray*}
This equality still holds after the change of the sign 
$a, n \to -a, -n$.
Therefore we have the equality \eqref{eq:recurrencef}. 

Next let us consider one of the next-to-boundary cases; 
$b=2, c=4, \tau(b)=1_{B}$.
Then we have 
\begin{eqnarray*}
& &  
\sum_{d=2, 4}q^{(m+1)w(d,2,4)}f_{a+2p'n, d, 2, m}(q) \\ 
& &  \qquad {}=
q^{(a/2+p'n)^{2}/t+m^2-(a/2+p'n)^{2}+(3m+2)/t+a/2+p'n+m} \\ 
& &  \qquad {}\times 
\Bigl\{
q^{(m+1)(y(4)-2\tau(4))-m}
\widetilde{S}_{m, 2-a/2-p'n}(q)+S_{m, 1-a/2-p'n}(q)
\Bigr\}
\end{eqnarray*}
Now the value $\tau(2)$ is equal to $1_{B}$, hence we have 
$\tau(4)=1_{B}$ or $2$, and then $y(4)-2\tau(4)=0$. 
Apply \eqref{eq:relS3} to the right hand side, and we get
\begin{eqnarray*}
& &  
\sum_{d=2, 4}q^{(m+1)w(d,2,4)}f_{a+2p'n, d, 2, m}(q) \\  
& &  {}=
q^{(a/2+p'n)^{2}/t+m^2-(a/2+p'n)^{2}+(3m+2)/t}
\Bigl\{
q^{a/2+p'n+m}S_{m+1, 1-a/2+p'n}(q)-S_{m,-a/2-p'n}(q)
\Bigr\}.
\end{eqnarray*}
Note that the power of $q$ in the coefficient of 
$S_{m, -a/2-p'n}(q)$ is invariant under the change of the sign 
$a, n \to -a, -n$. 
Therefore after the sum over $n \in \mathbb{Z}$ and $\epsilon=\pm$ 
in the right hand side of \eqref{eq:recurrencef},  
the function $S_{m, -a/2-p'n}(q)=S_{m, a/2+p'n}(q)$ disappears. 
Thus we obtain \eqref{eq:recurrencef}.   
\end{proof}

To finish the proof of Proposition \ref{PROP1},  
It suffices to show that 
\begin{eqnarray}
I_{r,a,b,m}^{(p,p')}(q)=
\sum_{d=b, \, b\pm2}
q^{m(\Delta_{r,d}-\Delta_{r,b}+\delta_{d,b})+\Delta_{r,b}-\Delta_{r,a}}
X_{a,d,b,m-1}(q). 
\label{eq:IandS}
\end{eqnarray}
{}From \eqref{eq:Xandf} and the recurrence relations 
\eqref{eq:relS1}--\eqref{eq:relS3}
we can check \eqref{eq:IandS} by direct calculation. 
\qed 

\noindent{\it Remark.}\quad
In \cite{FJMMT1}, 
we constructed a monomial basis of $M^{(p,p')}_{r,s}$ 
with $1< p'/p< 2$, using (2,1) primary field. 
To show that the monomials span the space, 
quadratic exchange relations were employed. 
In the process of rewriting the monomials, 
it was necessary to show the non-vanishing of 
a certain determinant. 
In the present case of (1,3) primary field, 
a similar set of quadratic relations can be written
explicitly.  
However it is not clear to us how to derive the spanning property
for the proposed set of monomials. 

{}From the next section, we will restrict to the 
case of unitary series $p'=p+1$. 

\section{Preliminaries on the fusion product}

\subsection{Fusion product}
In this section we fix our notation 
and collect the main properties of
the fusion product. 

Let $V_1,\ldots, V_n$ be cyclic representations of a 
Lie algebra $\g$ with cyclic vectors $v_1,\ldots,v_n$. 
Fix  $z=(z_1,\ldots,z_n)\in\C^n$ with $z_i\ne z_j$ for $i\ne j$.  Denote by $V_i(z_i)$ the evaluation representation of 
$\g\T\C[u]$, which is isomorphic to $V_i$ as vector space 
and the action is defined via the map 
$\g\T\C[u]\to\g$, $g\T u^j\mapsto z_i^j g$ ($g\in\g$). 
Recall (see \cite{FL}) that the fusion product
$V_1(z_1)*\ldots * V_n(z_n)$ is
the associated graded $\g\T \C[u]$-module 
with respect to the filtration $\{F_m\}_{m\ge 0}$ 
on the tensor product $V_1(z_1)\T \cdots \T V_n(z_n)$:
\begin{equation}
\label{fusfiltr}
F_m=\mathrm{span}\{\bigl(g_1\T u^{k_1}\cdots g_p\T u^{k_p}\bigr) 
(v_1\T\cdots\T v_n)\mid g_1,\cdots,g_p\in\g,~~k_1+\cdots +k_p\le m\}.
\end{equation}
If $[\mathfrak g,\mathfrak g]=\mathfrak g$, we have
\be
F_m=\mathfrak g\cdot\mathrm{span}\{\bigl(g_1\T u\cdots g_p\T u\bigr) 
(v_1\T\cdots\T v_n)\mid g_1,\cdots,g_p\in\g,~~p\le m\}.
\en
We set 
\begin{equation}
\label{degree}
(V_1(z_1)*\ldots * V_n(z_n))_m=F_m/F_{m-1},
\end{equation}
so that we have 
\be
V_1(z_1)*\ldots * V_n(z_n)=
\bigoplus_{m=0}^\infty (V_1(z_1)*\ldots * V_n(z_n))_m.
\en

The most important property of the fusion product is its independence on $z$ 
in some special cases (see \cite{FL, FF1, CL, FKL, FoL, Kedem, AK}).
Among such cases, we will need two cases: $\g=\slt$ and $V_i$ are
irreducible representations, and $\g=\slth'$ with $n=2$.
Note that, for an arbitrary Lie algebra $\mathfrak g$,
the fusion product of two representations $V_i$ $(i=1,2)$
with cyclic vectors $v_i$,
is independent of $z_1,z_2$ because
\bel
F_m=U(\mathfrak g)\cdot(G_m\otimes v_2)\hbox{ where }G_0=\C v_1,
\ G_{m+1}=G_m+\mathfrak g\cdot G_m.
\label{PBWFIL}
\enl
The filtration $G_m$ on $V_1$
is called the Poincar\'e-Birkhoff-Witt (PBW) filtration (see \ref{PBW}).

Consider the case where $\g=\slt$, $V_i=\pi_{a_i}$, 
and $v_i\in \pi_{a_i}$ is the highest weight vector. 
We write the corresponding fusion product as 
\begin{equation}
\label{fus}
\pi_{a_1} *\cdots *\pi_{a_n}
\end{equation}
as it is independent of the choice of $z$. 
The fusion product is also independent of the ordering
of the components $\pi_{a_i}$.
When $a_1=\cdots=a_n=a$, we use the shorthand notation 
$\pi_a^{*n}$.
The fusion product \eqref{fus} is a module over 
$\slt\T\C[u]$. 
We set $M^\alpha=\{v\in M\mid h_0v=\alpha v\}$ and 
\be
& & \ch_{q,z}\pi_{a_1} *\ldots *\pi_{a_n}
=\sum_{m,\alpha}q^mz^\alpha
\dim(\pi_{a_1} *\ldots *\pi_{a_n})^\alpha_m,
\\
& & 
\ch_{q}(\pi_{a_1} *\ldots *\pi_{a_n})^\alpha
=\sum_{m}q^m\dim(\pi_{a_1} *\ldots *\pi_{a_n})^\alpha_m.
\en
For example, 
\bea
& & {\rm ch}_{q,z}\pi_1^{* m}=
\sum_{-m\le l\le m\atop l\equiv m\bmod 2}
\qbin{m}{\frac{m+l}2}_{q}z^l,
\label{ch_pi1}\\
& & {\rm ch}_{q,z}\pi_2^{* m}=
\sum_{-m\le l\le m}S_{m,l}(q^{-1})z^{2l}. 
\label{ch_pi2}
\ena

\subsection{Exact sequences}
In this subsection we describe 
some exact sequences of fusion products of 
$\slt$-modules (see \cite{FF3}, Lemma $2.1$, Proposition $2.1$ and 
Remark $2.3$).

Let $A=(a_1,\cdots,a_n)$ be a sequence of non-negative integers 
with $a_1\le\ldots\le a_n$. 
Then for any $1\le i<  n$ there exists an 
$\slt\T\C[u]$-module $S_i(A)$ such that 
the following sequence is exact:
\begin{eqnarray*}
0
\longrightarrow
S_i(A)
\longrightarrow
\pi_{a_1}*\cdots *\pi_{a_n}
\longrightarrow
\pi_{a_1}*\cdots * 
\pi_{a_{i-1}}* \pi_{a_i-1} *\pi_{a_{i+1}+1}* 
\pi_{a_{i+2}}* \cdots * \pi_{a_n}
\longrightarrow
0\,.
\end{eqnarray*}
We will need two cases: $i=1$ and $i=n-1$.
First let $i=1$. Then
$$
S_1(A)\simeq \pi_{a_2-a_1}*\pi_{a_3}*\cdots *\pi_{a_n}.
$$
In order to have the above exact sequence to be degree preserving,
we must shift the $q$-degree of the highest weight vector of $S_1(A)$ 
(with respect to the operator $h\otimes u^0$) to $(n-1)a_1$. 
This gives the equality
\begin{eqnarray}
\ch_{q,z} \pi_{a_1}*\cdots *\pi_{a_n}
& =&  
q^{(n-1)a_1} \ch_{q,z} \pi_{a_2-a_1}*\pi_{a_3}*\cdots *\pi_{a_n}
\label{S_1}\\
& +& 
\ch_{q,z} \pi_{a_1-1}*\pi_{a_2+1}*\pi_{a_3}*
\cdots * \pi_{a_n}. 
\nonumber
\end{eqnarray}

Now let $i=n-1$. Then
$$S_{n-1}(A)\simeq \pi_{a_1} * \cdots * \pi_{a_{n-2}}\T \pi_{a_n-a_{n-1}}.$$
Therefore one has an exact sequence of $\slt\otimes\C[u]$-modules
\begin{multline}
\label{Sn}
0\to \pi_{a_1}*\cdots *\pi_{a_{n-2}}\T \pi_{a_n-a_{n-1}}
\to \pi_{a_1}*\cdots *\pi_{a_n}\to \\
\pi_{a_1}*\cdots * \pi_{a_{n-2}}* \pi_{a_{n-1}-1} *\pi_{a_n+1}\to 0
\end{multline}

We will also need an exact sequence 
involving different modules $S_i(A)$.
Namely let $A_i$ denotes the set of numbers
$$
(a_1, \ldots, a_{i-1}, a_i-1, a_{i+1}+1, a_{i+2}, \ldots,
a_n)
$$
reordered in the non-decreasing order. 
Then one has an exact sequence of $\slt\T\C[u]$-modules:
\begin{equation}
0\to \pi_{a_1}*\cdots *\pi_{a_{i-2}} *\pi_{a_{i-1}-a_i+a_{i+1}}* \pi_{a_{i+2}}*
\cdots *\pi_{a_n}\to S_i(A)\to S_i (A_{i-1})\to 0.
\end{equation}
For example, for $i=n-1$, we obtain
\begin{multline}
\label{Sind}
0\to \pi_{a_1}*\cdots *\pi_{a_{n-3}} *\pi_{a_{n-2}-a_{n-1}+a_n}\to
\pi_{a_1}* \cdots * \pi_{a_{n-2}}\T \pi_{a_n-a_{n-1}}\to\\
\pi_{a_1}* \cdots * \pi_{a_{n-2}-1}\T \pi_{a_n-a_{n-1}-1}\to 0.
\end{multline}

\subsection{$\slth$ and $\na_+$-homology}\label{HOMOLOGY}
We first settle our notation about $\slth$.  
Let 
$$
\slth=\slt\T\C[t,t^{-1}]\oplus \C K\oplus \C d,
$$
where $K$ is a central element and $[d,x_i]=-ix_i$, 
where we put $x_i=x\T t^i$ for $x\in\slt$.
Let 
$\na_+$ (resp. $\na_-$) be the 
nilpotent subalgebra of creation (resp. annihilation) 
operators generated by $f_0, e_{-1}$ (resp. $e_0$, $f_1$). 
For a positive integer $k$, 
we denote by $L_{l,k}$ ($0\le l\le k$) 
the set of integrable highest weight 
representations of $\slth$. 
We fix a highest weight vector $v_{l,k}\in L_{l,k}$.  
Then
$$
\na_- v_{l,k}=0, \ h_0v_{l,k}=l\,v_{l,k},\ 
K v_{l,k}=k\,v_{l,k},\ d\, v_{l,k}=0,\
L_{l,k}=U(\n_+)\cdot v_{l,k}.
$$
Representations $L_{l,k}$ are 
bi-graded by operators $d$ and $h_0$.
We set
$$
L_{l,k}=\bigoplus_{\al,n\in\Z} (L_{l,k})^\al_n, 
\quad 
(L_{l,k})^\al_n=\{v\mid dv=nv, h_0v=\al v\}.
$$
The Virasoro algebra acts on $L_{l,k}$ by the Sugawara operators:
$$L_n=\frac{1}{2(k+2)}\sum_{m\in\Z}:e_{n-m}f_m + f_{n-m} e_m +
\frac{1}{2} h_{n-m}h_m:,$$
where $:\ :$ is the normal ordering sign:
$$:x_i y_j:=
\begin{cases}
x_i y_j, \text{ if $i< j$;}\\
y_j x_i, \text{ if $i> j$;}\\
\frac{1}{2} (x_i y_i + y_i x_i), \text{ if $i=j$.}
\end{cases}
$$
The central charge is equal to $\frac{3k}{k+2}$. The conformal weight
$\triangle (l,k)$ of the highest weight vector $v_{l,k}$ is equal to
$\frac{l(l+2)}{4(k+2)}$: $L_0 v_{l,k}=\triangle(l,k) v_{l,k}$.

We now recall the homology result from \cite{FF4}.
For an $\slt\T\C[t]$-module $M$, 
we denote by $M'$ the $\slt\T\C[t^{-1}]$-module  
which is isomorphic to $M$ as a vector space and the
action  is defined via the isomorphism $x_i\mapsto x_{-i}$.
Let $a_1\le \ldots \le a_n\le k+1$.
Then (see Corollary 2.3 and Theorem 2.2 of [FF4]) for $0\le l\le k$ we have 
\begin{equation}
\label{vanish}
H_p (\na_+, (\pi_{a_1}*\cdots *\pi_{a_n})'\T L_{l,k})^0=0
\quad (p> 0), 
\end{equation}
where the superscript $0$ denotes 
the weight zero spaces with respect to the operator $h_0$.  
We also note that the characters of the 
zeroth homology groups of these spaces can be identified
with the Kostka polynomials (see \cite{FF4, FJKLM, SS}).

The following will be used later. 

\begin{lem}
For non-negative integers $j,k,l,m$ satisfying 
$0\le j\le k$, $0\le l\le k+1$,  we have 
\bea
& & H_p\bigl(\na_+,(\pi_1^{*m}\otimes \pi_j)'\otimes L_{l,k+1}
\bigr)^0=0
\quad (p> 0),
\label{vanish2}
\\
& & H_p\bigl(\na_+,(\pi_2^{*m}\otimes \pi_j)'\otimes L_{l,k+1}
\bigr)^0=0
\quad (p> 0).
\label{vanish3}
\ena
\end{lem}
\begin{proof}
First we show \eqref{vanish2}. We use a special case of (\ref{Sn}),
\be
0\to \pi_1^{*m}\otimes \pi_j
\to  \pi_1^{*(m+1)}* \pi_{j+1}
\to  \pi_1^{*m}* \pi_{j+2}
\to 0.
\en
Tensoring by $L_{l,k+1}$ we obtain 
\be
0\to (\pi_1^{*m}\otimes \pi_j)'\otimes L_{l,k+1}
\to  (\pi_1^{*(m+1)}* \pi_{j+1})'\otimes L_{l,k+1}
\to  (\pi_1^{*m}* \pi_{j+2})'\otimes L_{l,k+1}
\to 0.
\en
Assertion \eqref{vanish2} follows from 
the associated long exact sequence and 
\eqref{vanish}. 

To show \eqref{vanish3}, we first use a special case of (\ref{Sn}),
\be
0\to \pi_1*\pi_2^{*(m-1)}\otimes \pi_{j-1}
\to  \pi_1*\pi_2^{*m}* \pi_{j+1}
\to  \pi_1^{*2}*\pi_2^{*(m-1)}* \pi_{j+2}
\to 0.
\en
Arguing as above we obtain 
\be
H_p\bigl(\na_+,(\pi_1*\pi_2^{*(m-1)}\otimes \pi_{j-1})'
\otimes L_{l,k+1}
\bigr)^0=0
\quad (p> 0).
\en
Now applying a similar argument to a special case of (\ref{Sind}),
\be
0\to \pi_2^{*(m-1)}* \pi_{j+2}
\to  \pi_2^{*m}\otimes \pi_j
\to  \pi_1*\pi_2^{*(m-1)}\otimes \pi_{j-1}
\to 0,
\en
we obtain the desired result. 
\end{proof}

\section{The $\phi_{1,3}$ filtration.}

In what follows we deal with the Virasoro 
modules $M^{(p,p')}_{r,s}$ in the 
unitary case, $p=k+2, p'=k+3$.

Consider the decomposition 
of the tensor product of $\slth$-modules 
(the coset construction \cite{GKO})
\begin{equation}
\label{coset}
L_{i,1}\T L_{j,k}=\bigoplus_{0\le l\le k+1\atop 
l\equiv i+j\bmod 2}
M_{j+1,l+1}^{(k+2,k+3)}\T L_{l,k+1}.
\end{equation}
On the space $M_{j+1,l+1}^{(k+2,k+3)}$ we have 
two filtrations: 
the filtration defined by the (1,3) primary field,    
and the one induced from 
the fusion filtration on the left hand side. 
In this subsection, we show that these two filtrations 
coincide. 

For that purpose we use an operator 
identity due to \cite{L}, which we recall below. 
Consider the $\slth$ vertex operator associated with $\pi_2$, 
\bel
V^\sigma(z): L_{l,k+1}\T \pi_2\to
L_{l+\sigma,k+1}\T z^{\triangle(l+\sigma,k+1)-\triangle(l,k+1)-
\triangle(2,k+1)}\C((z)).
\label{INTERTWINING}
\enl
Here $\sigma=-2,0,2$  
and $\C((z))$ is the space of Laurent series.
We fix a weight basis $v_\tau\in \pi_2$ 
with $h v_\tau=\tau v_\tau$ ($\tau=-2,0,2$), 
and write
$$
V^\sigma(z)=(V^\sigma_{-2}(z), V^\sigma_0(z), V^\sigma_2(z)),
$$
where 
$V^\sigma_{\tau}(z)(u)=V^\sigma(z)(u\T v_\tau).$
We have the Fourier expansion 
\begin{eqnarray*}
& & V^\sigma_\tau(z)
=\sum_{n\in \Z+\triangle(l,k+1)-\triangle(l+\sigma,k+1)}
z^{-n-\triangle(2,k+1)} V^\sigma_{\tau,n},
\\
& & V^\sigma_{\tau,n}: (L_{l,k+1})^\al_\beta\to 
(L_{l+\sigma,k+1})^{\al+\tau}_{\beta-n}.
\end{eqnarray*}
We will also need the (1,3) field for the Virasoro modules.  
Abbreviating 
$\phi^{(s+\sigma,s)}(z)$ to  $\phi^\sigma(z)$ we write 
\begin{eqnarray*}
& & \phi^\sigma(z)=
\sum_{n\in \Z+\triangle_{r,s}-\triangle_{r,s+\sigma}}
\phi^\sigma_{n} z^{-n-\triangle_{1,3}},\\
& & \phi^\sigma_{n}: (M^{(p,p')}_{r,s})_\beta\to
(M^{(p,p')}_{r,s+\sigma})_{\beta-n}.
\end{eqnarray*}
In what follows 
we suppose that some normalization of vertex operators
$V^\sigma_\tau(z)$ and $\phi^\sigma(z)$ is fixed.

For $x\in\slt$, set $x(z)=\sum_{n\in\Z} x_n z^{-n-1}$. 
Introduce further the current
$$
j(z)=(j_{-2}(z), j_0(z),j_2(z))
$$ 
acting on the tensor product $L_{i,1}\T L_{j,k}$ by 
\begin{eqnarray*}
& & j_{-2}(z)=k f(z)\T \Id-\Id\T f(z), \\
& & j_0(z)=k h(z)\T \Id-\Id\T h(z), \\
& & j_2(z)=k e(z)\T \Id-\Id\T e(z).  
\end{eqnarray*}

The following proposition is proved in \cite{L}.
\begin{prop}
\label{lash}
There exist non-vanishing constants $c_{-2}, c_0, c_2$ 
such that, 
with respect to the identification \eqref{coset}, 
the following equality holds:
\begin{equation}
\label{psi}
j_\tau(z)=c_{-2}\phi^{-2}(z)\T V^{-2}_\tau(z)+ c_0\phi^0(z)\T V^0_\tau(z)
+c_2\phi^2(z)\T V^2_\tau(z).
\end{equation}
\end{prop}

Noting that $\Delta_{1,3}+\Delta(2,k+1)=1$, 
we set 
$$
j_\tau^\sigma(z)=
c_\sigma\phi^\sigma(z)\T V_\tau^\sigma(z)
=\sum j^\sigma_{\tau,n}z^{-n-1}, 
$$
so that 
$j_\tau(z)=j_\tau^{-2}(z)+j_\tau^0(z)+j_\tau^2(z)$. 

\begin{lem}
\label{homog}
For each $\gamma\in\Z_{\geq0}$, we have the equality
\be
& & \U(\slth)\cdot \mathrm{span} 
\{j_{\tau_1,i_1}\cdots j_{\tau_\gamma,i_\gamma}
(|r,s\rangle\T v_{r-1,k+1})
\mid \tau_\al=-2,0,2, i_\al\in\Z\}\\
& & =\U(\slth)\cdot\mathrm{span}\{j^{\sigma_1}_{\tau_1,i_1}\cdots 
j^{\sigma_\gamma}_{\tau_\gamma,i_\gamma}(|r,s\rangle\T v_{r-1,k+1})
\mid\tau_\al,\sigma_\beta=-2,0,2, i_\al\in\Z\}.
\en
\end{lem}
\begin{proof}
We note that for any
set of vectors $\{u_l\in L_{l,k+1}\}_{l=0}^{k+1}$ 
there exist elements
$\{x_l\in U(\slth)\}_{l=0}^{k+1}$ such
that  
$$
x_l\cdot u_{l'}=\delta_{l,l'} u_{l'}.
$$
Using this fact and the intertwining property of the vertex operator
\eqref{INTERTWINING}, Lemma can be proved by induction on $\gamma$.
\end{proof}

\begin{prop}
\label{mprop}
Fix $r,s$ and $k$. For $\sigma_1,\ldots,\sigma_\gamma=-2,0,2$, set
\be
& & W_{\sigma_1,\ldots,\sigma_\gamma}=U(\widehat{\mathfrak{sl}}_2)\cdot\\
& & \times\hbox{\rm span}
\{j^{\sigma_1}_{\tau_1,i_1}\cdots j^{\sigma_\gamma}_{\tau_\gamma,i_\gamma}
(|r,s\rangle\otimes v_{s-1,k+1})|\tau_\alpha=-2,0,2,i_\alpha\in\Z\}
\en
and
\be
& & \widetilde W_{\sigma_1,\ldots,\sigma_\gamma}=
\hbox{\rm span}
\{\phi^{\sigma_1}_{i_1}\cdots \phi^{\sigma_\gamma}_{i_\gamma}
|r,s\rangle|i_\alpha\in\Z\}\otimes L_{s+\sigma_1+\cdots+\sigma_\gamma-1,k+1}.
\en
Then, the equality $W_{\sigma_1,\ldots,\sigma_\gamma}
=\widetilde W_{\sigma_1,\ldots,\sigma_\gamma}$ holds.
\end{prop}
\begin{proof}
Let $a=s+\sigma_1+\cdots+\sigma_\gamma$.
We first note that
\be
W_{\sigma_1,\ldots,\sigma_\gamma}\subset
\widetilde W_{\sigma_1,\ldots,\sigma_\gamma}
\subset
M^{(k+2,k+3)}_{r,a}\otimes L_{a-1,k+1}.
\en
We show the equality
\be
W_{\sigma_1,\ldots,\sigma_\gamma}\cap
\left((M^{(k+2,k+3)}_{r,a})_n\otimes L_{a-1,k+1}\right)
=\widetilde W_{\sigma_1,\ldots,\sigma_\gamma}\cap
\left((M^{(k+2,k+3)}_{r,a})_n\otimes L_{a-1,k+1}\right).
\en
Note that the equality follows if we show that
for any $n_\alpha$ such that $n_1+\cdots+n_\gamma=n-\Delta_{r,s}$,
\be
\phi^{\sigma_1}_{-n_1}\cdots\phi^{\sigma_\gamma}_{-n_\gamma}|r,s\rangle
\otimes v_{a-1,k+1}\in W_{\sigma_1,\ldots,\sigma_\gamma}
\en
holds. We prove this statement by induction on $n$.

Before starting, we note that
in view of the equality $\Delta_{1,3}+\Delta(2,k+1)=1$, we obtain
\be
j^{\sigma_\alpha}_{\tau_\alpha,i_\alpha-n_\alpha}=
\sum_{i-n=i_\alpha-n_\alpha}
\phi^{\sigma_\alpha}_{-n}\otimes V^{\sigma_\alpha}_{\tau_\alpha,i}.
\en

For $n=\Delta_{r,a}$ the statement is equivalent to
\be
|r,a\rangle\otimes v_{a-1,k+1}\in W_{\sigma_1,\ldots,\sigma_\gamma}.
\en
This is a consequence of the fact that
for any $s_1,s_2$ such that $\sigma=s_1-s_2=-2,0,2$
one can find unique $\tau,i,i_1,i_2$ $(i_1+i_2=i)$ such that
\be
|r,s_1\rangle\otimes v_{s_1-1,k+1}\in
\C j^{\sigma}_{\tau,i}(|r,s_2\rangle\otimes v_{s_2-1,k+1})=
\C \phi^\sigma_{i_1}|r,s_2\rangle
\otimes V^{\sigma}_{\tau,i_2}v_{s_2-1,k+1}.
\en


Now suppose the statement is proved for all $n\leq n_0$. Fix
$\tau_\alpha,i_\alpha$ such that
\be
v_{a-1,k+1}\in
\C V^{\sigma_1}_{\tau_1,i_1}\cdots
V^{\sigma_\gamma}_{\tau_\gamma,i_\gamma}v_{s-1,k+1}.
\en
Choose $n_\alpha$ such that $n_1+\cdots+n_\gamma=n_0+1-\Delta_{r,s}$
and consider
\be
w_{n_1,\ldots,n_\gamma}=j^{\sigma_1}_{\tau_1,i_1-n_1}\cdots
j^{\sigma_\gamma}_{\tau_\gamma,i_1-n_\gamma}
(|r,s\rangle\otimes v_{s-1,k+1}).
\en
We have
\be
w_{n_1,\ldots,n_\gamma}=
c\phi^{\sigma_1}_{-n_1}\cdots\phi^{\sigma_\gamma}_{-n_\gamma}
|r,s\rangle\otimes v_{a-1,k+1}+\sum_{n=\Delta_{r,a}}^{n_0}
\tilde w_n,
\en
where $c$ is a non-zero constant and
\be
\tilde w_n\in
\widetilde W_{\sigma_1,\ldots,\sigma_\gamma}\cap
\left((M^{(k+2,k+3)}_{r,a})_n\otimes L_{a-1,k+1}\right).
\en
The statement follows from the induction hypothesis.
\end{proof}

We now define the $\phi_{1,3}$ filtration on $M^{(k+2,k+3)}_{r,s}$.
Note that when $p'=p+1$, the integer $b(r,s)$ determined
by the condition (2.6) reads
\be
b(r,s)=\begin{cases}r& \hbox{ if $r-s$ is even:}\\
r+1& \hbox{ otherwise.}
\end{cases}
\en
Now, consider the filtration
\bel
E_m(M^{(k+2,k+3)}_{r,s})=
{\rm span}\{\phi^{\sigma_1}_{i_1}\cdots
\phi^{\sigma_\gamma}_{i_\gamma}|r,b(r,s)\rangle
|\gamma\leq m,\sigma_\alpha=-2,0,2\}\cap M^{(k+2,k+3)}_{r,s}.
\label{13filt}
\enl

Let
\be
{\rm Gr}^E_m(M^{(k+2,k+3)}_{r,s})
=E_m(M^{(k+2,k+3)}_{r,s})/E_{m-1}(M^{(k+2,k+3)}_{r,s}).
\en

Recall the decomposition (see (\ref{degree}))
$$
L_{i,1}* L_{r-1,k}=\bigoplus_{m=0}^\infty 
(L_{i,1}* L_{r-1,k})_m.
$$
\begin{prop}
\label{fusphi}
Under the identification $(\ref{coset})$, 
we have an isomorphism of $\slth$-modules
$$
(L_{i,1}* L_{r-1,k})_m=
\bigoplus_{1\le s\le k+2\atop s\equiv r+i \bmod 2}
{\rm Gr}^E_m(M^{(k+2,k+3)}_{r,s})\T L_{s-1,k+1}.
$$
\end{prop}
\begin{proof}
We use the definition of the fusion product 
\eqref{fusfiltr}
of $\slth$-modules, choosing 
$$
V_1=L_{i,1},~~V_2=L_{r-1,k}, ~~
v_1=v_{i,1},~~ v_2=v_{r-1,k},~~
z_1=k,~~z_2=-1.
$$
Then the fusion filtration $F_m$ 
on this tensor product is given by
$$
F_m=\U(\slth)\cdot \mathrm{span}\{
j_{\tau_1,i_1}\cdots j_{\tau_\gamma,i_\gamma}
(|r,r+i\rangle\T v_{s-1,k+1})
\mid \tau_\al=-2,0,2, ~~i_\al\in\Z,~~
\gamma\le m\}. 
$$
Now our proposition follows from Lemma \ref{homog} and
Proposition \ref{mprop}.
\end{proof}

\section{Decomposition of $\slth$ fusion products}

Our goal in this section is to determine the character of 
${\rm Gr}^E_m(M^{(k+2,k+3)}_{r,a})$ 
which appear in the decomposition of 
the fusion product  $L_{i,1} * L_{j,k}$. 
We show that it is given in terms of
$I^{(k+2,k+3)}_{a,b,r,m}(q)$ introduced in \eqref{eq:Im}. 

\subsection{The functor $I_k$.}\label{FUNCTOR}
Set
\be
\mathfrak{g}_\pm=\slt\otimes\C[t^{\mp1}]\oplus\C d,
\quad
\g=\slth. 
\en
Let $V$ be a $\g_-$-module with $\Z_{\ge 0}\times\Z$-grading, 
\be
V=\oplus_{s\ge 0}\oplus_{\alpha\in\Z}V^{\alpha}_s,
\quad V^{\alpha}_s=\{v\in V\mid d\,v=s\,v,~~h_0v=\alpha\,v\}.
\en
The induced $\g$-module 
\be
\mathrm{Ind}_{\g_-}^{\g}V
=U(\g)\otimes_{U(\g_-)}V
\en
is also $\Z_{\ge 0}\times \Z$-graded. 
We will consider its maximal integrable quotient.

Let $k$ be a positive integer.  
Since the $q$-degree is bounded from below, 
the operator $e(z)^{k+1}=\sum_ne^{(k+1)}_nz^{-n-1}$
has a well-defined action on ${\rm Ind}^{\mathfrak g}_{\mathfrak g_-}V$.
Let $\mathcal I_k$ be the $\mathfrak g$-submodule of
${\rm Ind}^{\mathfrak g}_{\mathfrak g_-}V$
generated by
\be
(K-k)\cdot{\rm Ind}^{\mathfrak g}_{\mathfrak g_-}V+\sum_n
e^{(k+1)}_n\cdot{\rm Ind}^{\mathfrak g}_{\mathfrak g_-}V.
\en
We define
\be
I_k(V)={\rm Ind}^{\mathfrak g}_{\mathfrak g_-}V/\mathcal I_k.
\en
Then $I_k(V)$ is an integrable $\g$-module of level $k$. 
Moreover, any homomorphism $V\to L$ of $\g_-$-modules to 
an integrable $\g$-module $L$ of level $k$ 
extends to a homomorphism of $\g$-modules $I_k(V)\to L$. 

Since $I_k(V)$ is integrable of positive level $k$, 
it has the decomposition 
\be
I_k(V)=\bigoplus_{l=0}^k I_k^l(V)\otimes L_{l,k},
\en
where
\bel
I_k^l(V)=\mathrm{Hom}_{\g'}(L_{l,k},I_k(V))
\label{INDUC}
\enl
is the space of highest weight vectors of weight $l$. 
Here we have set 
$\g'=[\g,\g]=\slt\otimes\C[t,t^{-1}]\oplus\C K$. 
The space $I_k^l(V)$ carries a grading by $d$. 
In the next subsection we compute the character of $I_k^l(V)$ 
for some $V$. 

\subsection{Tensor product as induced module}

In this subsection, we show that the tensor product module
$L_{0,1}\otimes L_{j,k}$
can be realized as  $I_{k+1}(L_{0,1}\otimes \pi_j)$
(see Proposition \ref{prop:LLIk} below). 
We first prepare two Lemmas.  

Let $\omega$ denote the involutive automorphism of $\slth$ 
given by
\be
e_i\mapsto f_{-i},~~
f_i\mapsto e_{-i},~~
h_i\mapsto -h_{-i},~~
K\mapsto -K,~~
d\mapsto -d. 
\en
For a $\g_-$-module $V$ defined by 
$\rho_-:\g_-\to\mathrm{End}(V)$, 
let $V^{\omega}$ denote the 
$\g_+$-module structure on $V$ given by
$\rho_+=\rho_-\circ\omega:\g_+\to\mathrm{End}(V)$. 

We consider the situation where $V$ is the fusion product
of $\mathfrak{sl}_2$-modules or their tensor products. We identify
$u$ with $t$ and define the degree operator $d$ appropriately
(see Proof of Proposition 5.4). In this case, $V^\omega$ and $V'$ (see 3.3)
are isomorphic as $\mathfrak{sl}_2\otimes\C[t^{-1}]$-module.

Recall the definition of $\mathfrak n_+$ given in \ref{HOMOLOGY}.

\begin{lem}\label{H_0}
Notation being as above, we have 
\begin{eqnarray*}
\ch_q I_k^l(V)
& =& 
\ch_{q^{-1}} H_0(\na_+, V^{\omega}\otimes L_{l,k})^0. 
\end{eqnarray*}
\end{lem}
\begin{proof}
In order to use the reciprocity law, we rewrite (\ref{INDUC}) as
\be
I_k^l(V)\simeq
\mathrm{Hom}_{\g'} (I_k(V), L_{l,k})^*, 
\en
where $*$ denotes the restricted dual. 
Setting $\g'_\pm=\slt\otimes\C[t^{\mp 1}]$, we have 
\be
\mathrm{Hom}_{\g'} (I_k(V),L_{l,k})^*
& \simeq & 
\mathrm{Hom}_{\g'} 
(\mathrm{Ind}_{\g'_-}^{\g'}V, L_{l,k})^*
\\
& \simeq & 
\mathrm{Hom}_{\g'_-}(V, L_{l,k})^*
\\
& \simeq & 
\Bigl((L_{l,k}\otimes V^*)^{\g_-'}\Bigr)^*
\\
& \simeq & 
(V\otimes L_{l,k}^*)/\g_-'(V\otimes L_{l,k}^*)
\\
& \simeq & 
H_0(\g_-',V\otimes L_{l,k}^*).
\en
Note that as $\g_+'$-modules we have 
$(L_{l,k}^*)^\omega\simeq L_{l,k}$. 
Therefore, we obtain
\be
{\rm ch}_qH_0(\mathfrak g'_-,V\otimes L^*_{l,k})
& \simeq& {\rm ch}_{q^{-1}}H_0(\mathfrak g'_+,V^\omega\otimes L_{l,k})\\
& \simeq& {\rm ch}_{q^{-1}}H_0(\mathfrak n_+,V^\omega\otimes L_{l,k})^0.
\en
Now our lemma follows. 
\end{proof}


\begin{lem}\label{Euler}
Let $V$ be a bi-graded $\g_+$-module. Then
\bea
& & \sum_{p=0}^\infty (-1)^p \ch_q H_p(\na_+, V\T L_{l,k})^0
\label{char}\\
& & 
\quad=
\sum_{\la\in\Z} 
q^{(k+2)\la^2-(l+1)\la}
\left(\ch_q V^{-2\la(k+2)+l}-\ch_q V^{-2\la(k+2)+l+2}\right).
\nonumber
\ena
\end{lem}
\begin{proof}
We first prepare our notation concerning 
the Weyl group $W$ of $\slth$.
It is generated by simple reflections $s_0,s_1$.
The length of $w\in W$ is denoted by $\ell(w)$. 
Let $\mathfrak h=\C h_0\oplus \C K\oplus\C d$ be the Cartan subalgebra.
Define $(i,k,m)\in\mathfrak h^*$ by $(i,k,m)h_0=i,
(i,k,m)K=k,(i,k,m)d=m$. We have
\be
& & s_1(i,k,m)=(i,k,m)-i\alpha_1,\quad\alpha_1=(2,0,0),\\
& & s_0(i,k,m)=(i,k,m)-(k-i)\alpha_0,\quad\alpha_0=(-2,0,-1).
\en
Let $\rho=(1,2,0)$. We write
$w*\alpha=w(\alpha+\rho)-\rho$ for the shifted action of the Weyl group
on $\mathfrak h^*$.

%
Consider the BGG-resolution of $L_{l,k}$ (see \cite{BGG, Kum})
\bea
& & 0
\longleftarrow
L_{l,k}
\longleftarrow
M_0
\longleftarrow
M_1
\longleftarrow
\cdots,
\label{fBGG}
\\
& & 
M_p=\bigoplus_{\ell(w)=p} M(w*(l,k,0)),
\nonumber
\ena
where $M(\mu)$ is the Verma module with highest weight
$\mu$. 
We recall that, for a general Lie algebra $\mathfrak{a}$,    
a $U(\mathfrak{a})$-module $V$ 
and a free $U(\mathfrak{a})$-module $M$, 
the $U(\mathfrak{a})$-module 
$V\otimes M$ is isomorphic to 
$V^{\mathrm{triv}}\otimes M$ where 
$V^{\mathrm{triv}}$ is $V$ with trivial $\mathfrak{a}$-action. 
Hence, tensoring (\ref{fBGG}) by $V$, 
we obtain a $U(\mathfrak n_+)$-free resolution for $V\T L_{l,k}$.
Therefore, the following complex counts 
$H_*(\na_+,V\T L_{l,k})^0$:
\begin{equation}
\label{Fp}
0\longleftarrow (\C\T_{U(\na_+)} (V\otimes M_0))^0\longleftarrow
(\C\T_{U(\na_+)} (V\otimes M_1))^0\longleftarrow
\cdots
\end{equation}
We can rewrite
(\ref{Fp}) as
\begin{equation*}
0\longleftarrow ((\C\T_{U(\na_+)} M_0)\T V)^0\longleftarrow
((\C\T_{U(\na_+)} M_1)\T V)^0\longleftarrow
\cdots
\end{equation*}
By the Euler-Poincar{\'e} principle, 
\be
\sum_{p=0}^\infty(-1)^p\ch_qH_p(\na_+,V\otimes L_{l,k})^0
& =& 
\sum_{p=0}^\infty(-1)^p
\ch_q\bigl((\C\otimes_{U(\na_+)}M_p)\otimes V\bigr)^0
\\
& =& 
\sum_{w\in W}(-1)^{\ell(w)}
q^{(w*(l,k,0))(d)}
\ch_q V^{-(w*(l,k,0))(h_0)}.
\en

Lemma follows by using the following formulas for the 
shifted action where $\lambda\in\Z$ ,
\begin{gather}
\label{W}
(s_1s_0)^\lambda * (l,k,m)
=(-2(k+2)\lambda+l,k,m+(k+2)\lambda^2-(l+1)\lambda),
\\  \nonumber
s_0(s_1s_0)^{\lambda-1} * (l,k,m)
=(-l-2+2(k+2)\lambda, k, m+(k+2)\lambda^2-(l+1)\lambda).
\end{gather}
\end{proof}

\begin{cor}\label{Euler2}
If $H_p(\na_+,V^{\omega}\otimes L_{l,k})=0$ holds for $p> 0$, 
then
\be
\ch_qI^l_k(V)=
\sum_{\lambda\in\Z}q^{-(k+2)\lambda^2+(l+1)\lambda}
\left(\ch_qV^{2(k+2)\lambda-l}-\ch_qV^{2(k+2)\lambda-l-2}\right). 
\en
\end{cor}
\begin{proof}
This follows from Lemma \ref{H_0} and Lemma \ref{Euler},
with $q,V$ replaced by $q^{-1},V^\omega$ noting that 
$\ch_{q^{-1}}(V^{\omega})^\alpha=\ch_q V^{-\alpha}$. 
\end{proof}

\begin{prop}\label{prop:LLIk}
We have an isomorphism of $\slth$-modules 
\be
L_{0,1}\otimes L_{j,k}\simeq I_{k+1}(L_{0,1}\otimes \pi_j), 
\en
where we regard $\pi_j$ as an  
$\slt\otimes\C[t]\oplus\C d$-module by letting 
$\slt\otimes t\C[t]\oplus\C d$ act as $0$. 
\end{prop} 
\begin{proof}
As we noted in \ref{FUNCTOR}, there  exists a natural surjective homomorphism 
\be
I_{k+1}(L_{0,1}\otimes\pi_j)\to L_{0,1}\otimes L_{j,k}. 
\en
Therefore, it suffices to check that the multiplicities of 
irreducible representations in the decomposition
of left and right hand sides coincide.

Recall that $L_{0,1}$ is an inductive limit of its 
Demazure submodules, which 
are isomorphic to fusion products 
of $2$-dimensional representations \cite{FF2}. 
Introduce the action of $d$ on $\pi_1^{*2N}$ by setting 
\be
d\,v=(N^2-m)\,v, \quad v\in (\pi_1^{*2N})_m,
\en
and denote by $\widetilde{\pi_1^{*2N}}$ the resulting 
$\slt\otimes\C[t]\oplus\C d$-module. 
Then we have an isomorphism of 
$\slt\otimes\C[t]\oplus\C d$-modules   
\bea
L_{0,1}=\lim_{N\to\infty} \widetilde{\pi_1^{*2N}}.
\label{Demaz}
\ena
Since inductive limit commutes with tensor product
and is an exact functor, we have
\be
I_{k+1}(L_{0,1}\otimes\pi_j)
=\lim_{N\to\infty}
I_{k+1}(\widetilde{\pi_1^{*2N}}\otimes\pi_j). 
\en
Therefore 
\begin{eqnarray*}
\ch_q I_{k+1}(L_{0,1}\otimes \pi_j)
& =& \lim_{N\to\infty}
\ch_q I_{k+1}(\widetilde{\pi_1^{*2N}}\otimes\pi_j). 
\end{eqnarray*}

In view of the vanishing of homology \eqref{vanish2}, 
we can apply Corollary \ref{Euler2}. 
We obtain 
\be
& & \ch_q I_{k+1}^l(\widetilde{\pi_1^{*2N}}\otimes\pi_j). 
\\
& & \quad
=
\sum_{\la\in\Z} q^{-(k+3)\la^2+(l+1)\la}
\left(
\ch_q (\widetilde{\pi_1^{*2N}}
\otimes \pi_j)^{2(k+3)\lambda-l}-\ch_q (\widetilde{\pi_1^{*2N}}
\otimes \pi_j)^{2(k+3)\lambda-l-2}
\right)
\\
& & \quad
=
\sum_{\la\in\Z} q^{-(k+3)\la^2+(l+1)\la}
\left(
\ch_q (\widetilde{\pi_1^{*2N}})^{2(k+3)\lambda-l+j}-
\ch_q (\widetilde{\pi_1^{*2N}})^{2(k+3)\lambda-l-2-j}
\right).
\en
On the other hand, we have from \eqref{ch_pi1}
$$
\ch_q(\widetilde{\pi_1^{*2N}})^\al
=
q^{N^2}{\bininv{2N}{\frac{2N+\al}{2}}}
=q^{\alpha^2/4}{\bin{2N}{\frac{2N+\al}{2}}}.
$$
Hence we find 
\be
& & 
\ch_q I_{k+1}^l(\widetilde{\pi_1^{*2N}}
\otimes\pi_j)
\\
& & \quad
=
q^{\frac{(l-j)^2}{4}}\Biggl(
\sum_{\la\in\Z}
q^{(k+2)(k+3)\la^2+((k+3)(j+1)-(k+2)(l+1))\la}
\bin{2N}{\frac{2N-l+j}{2}+(k+3)\la}
\\
& & \quad-
\sum_{\la\in\Z} 
q^{(k+2)(k+3)\la^2-((k+3)(j+1)+(k+2)(l+1))\la+(j+1)(l+1)}
\bin{2N}{\frac{2N-l-j-2}{2}+(k+3)\la}\Biggr).
\en
The last expression in the above formula 
coincides with the ABF finitization 
of the Virasoro minimal model 
$M^{(k+2,k+3)}_{j+1,l+1}$ up to a power
$q^{\triangle_{j+1,l+1}-(l-j)^2/4}$ (see \cite{ABF}).
This power comes from conformal dimensions of
highest weight vectors of $\slth$-modules:
$$
\triangle_{j+1,l+1}-\frac{(l-j)^2}{4}=
\triangle(j,k)-\triangle(l,k+1),
$$
the right hand side being 
the difference between $L_0$-grading and
$d$-grading. 
Letting $N\to\infty$ we conclude that
\be
\ch_qI_{k+1}^l(L_{0,1}\otimes\pi_j)
=q^{\Delta(l,k+1)-\Delta(j,k)}\times
\ch_qM^{(k+2,k+3)}_{j+1,l+1}. 
\en
This completes the proof of Proposition. 
\end{proof}

\subsection{PBW-filtration on $L_{0,1}$}\label{PBW}

Consider the fusion filtration $F_m$ on the tensor product
$L_{0,1}\otimes L_{j,k}$ given by (\ref{PBWFIL}).
In order to study $F_m/F_{m-1}$, 
we consider in this subsection the structure of the 
PBW-filtration $G_m$ on $L_{0,1}$. 
Denote the associated graded space by $Q_m=G_m/G_{m-1}$. 
Then $Q_m$ inherits the $q,z$-grading from $L_{0,1}$.
Note that for any $n\geq0$
\be
(Q_m)^{2m}_n=(G_m)^{2m}_n.
\en
It is known \cite{FS} that this space has a monomial basis
\bel
e_{-n_1}\cdots e_{-n_m}v_{0,1}
\label{MONOMIAL}
\enl
where $n_j\geq n_{j+1}+2$, $n_m> 0$ and $\sum_{j=1}^mn_j=n$.
In particular, we have
\be
\ch_q(Q_m)^{2m}=\frac{q^{m^2}}{(q)_m}.
\en
There is a canonical action of $\mathfrak g_-$ on $Q_m$. 
Define a $\mathfrak g_-$-invariant filtration on $Q_m$ by setting
\be
J_n(Q_m)=\sum_{i\leq n}U(\mathfrak g_-)(Q_m)^{2m}_{m^2+i}.
\en
\begin{prop}\label{Qm}
We have
\be
Q_m=\cup_{n=0}^\infty J_n(Q_m).
\en
The associated graded space ${\rm Gr}^J_n(Q_m)$
is isomorphic to a direct sum of the fusion product $\pi_2^{*m}$.
Each monomial vector \eqref{MONOMIAL} generates an
$\mathfrak{sl}_2\otimes\C[t]$-module isomorphic to $\pi_2^{*m}$.
\end{prop}
Obviously, $\cup_{n=0}^\infty J_n(Q_m)\hk Q_m$. 
For the proof of this proposition we need the following lemma.

\begin{lem}\label{CHARACTER}
We have
\be
\ch_q Q_m=\frac{q^{m^2}}{(q)_m}\ch_{q^{-1}} \pi_2^{*m}.
\en
\end{lem}
The proof of this lemma is given below.

We also need the dual functional realization for the space $Q_m$.
Set $K_m=\C[z_1,\ldots,z_m]$.
Consider the restricted dual space
\be
Q_m^*=\oplus_n(Q_m)_n^*.
\en
For $\varphi\in Q_m^*$ one can define an element
$F_\varphi\in K_m\otimes(\mathfrak{sl}_2^{\otimes m})^*$ by
\bel
F_\varphi(x^{(1)}\otimes\cdots\otimes x^{(m)})
=\varphi(x^{(1)}_+(z_1)\cdots x^{(m)}_+(z_m)v_{0,1}).
\enl
Here for $x\in\mathfrak{sl}_2$ we set $x_+(z)=\sum_{n\geq1}x_{-n}z^n$.
The map $\varphi\mapsto F_\varphi$ is injective. We will give a
characterization of the image of this mapping. First we recall some results
of \cite{FF1} on the dual functional realization of the fusion product.

Note the $\mathfrak{sl}_2$ decomposition
\be
\pi_2\otimes\pi_2=\pi_0\oplus\pi_2\oplus\pi_4.
\en
Define
\bel
V_m& =& \{F(z_1,\ldots,z_m)\in K_m\otimes \pi_2^{\otimes m}|\nonumber\\
& & \hskip5pt
F|_{z_i=z_j}\in\sigma^{(i,j)}\left(\pi_0\otimes\pi_2^{\otimes(m-2)}\right),
\label{V1}\\
& & \hskip5pt
(\partial_{z_i}F)|_{z_i=z_j}
\in\sigma^{(i,j)}\left((\pi_0\oplus\pi_2)\otimes\pi_2^{\otimes(m-2)}\right)\}.
\label{V2}
\enl
Here
\be
\sigma^{(i,j)}(v_i\otimes v_j\otimes v_1\otimes\cdots\otimes\widehat{v_i}
\otimes\cdots\otimes\widehat{v_j}\otimes\cdots\otimes v_n)=
v_1\otimes\cdots\otimes v_n.
\en
The space $V_m$ is a $\mathfrak g_-$-module by the action
\be
x_n\cdot F=\sum_{j=1}^mz_j^n\varpi_j(x)F.
\en
Here $\varpi_j(x)$ is the action of $x\in\mathfrak{sl}_2$
on the $j$-th component of the tensor product $\pi_2^{\otimes m}$.
The homogeneous degree of elements in $K_m$ is counted by $-d\in\mathfrak g_-$.
We denote by $(K_m)_n$ the subspace with homogeneous degree $n$.

Multiplication by elements of $K_m$ commutes with this action.
Let $K_m^0$ be the maximal ideal of $K_m$ generated by $z_1,\ldots,z_m$.
The following result is proved in \cite{FF1} (see Proposition 4.1
and Theorem 4.1).

\begin{prop}\label{FUSIONDUAL}
We have an isomorphism of $\mathfrak{sl}_2\otimes\C[t]$-modules
\be
V_m/K_m^0V_m\simeq (\pi_2^{*m})^*.
\en
\end{prop}
We need a variant of this result in the symmetric case.
The symmetric group $\mathfrak S_m$ acts on $V_m$ by
\bel
(\sigma F)(z_1,\ldots,z_m)=
(1\otimes\sigma)\cdot F(z_{\sigma(1)},\ldots,z_{\sigma(m)}),
\label{ACTION}
\enl
where the action of $\sigma$ on $\pi_2^{\otimes m}$ is
$\sigma(v_1\otimes\cdots\otimes v_m)=
v_{\sigma^{-1}(1)}\otimes\cdots\otimes v_{\sigma^{-1}(m)}$.
Set
\be
W_m=V_m\cap(K_m\otimes\pi_2^{\otimes m})^{\mathfrak S_m}.
\en
Let $S_m=K_m^{\mathfrak S_m}$ and $(S_m)_n=(K_m)_n^{\mathfrak S_m}$.
Note that
\bea
(W_m)^{-2m}=\prod_{i< j}(z_i-z_j)^2S_m\otimes v^{\otimes m}, 
\label{eq:degree-correspondence-W}
\ena
where $v$ is the lowest weight vector of $\pi_{2}$.

Set
\be
S_m^0=K_m^0\cap S_m.
\en
This is the ideal of $S_m$ consisting of symmetric polynomials vanishing
at the origin.
\begin{lem}\label{NATURAL}
The natural map
\be
W_m/S_m^0W_m\rightarrow V_m/K_m^0V_m
\en
is an isomorphism.
\end{lem}
\begin{proof}
First we prove 
\be
W_m\cap(K_m^0V_m)\subset S_m^0W_m.
\en
It is enough to show that for any $n> 0$
\bel
{\rm Sym}(gF)\in S_m^0W_m\quad\hbox{for any $g\in(K_m^0)_n$ and $F\in V_m$.}
\label{SYMMETRIZATION}
\enl
Here ${\rm Sym}\,F=\frac1{m!}\sum_{\sigma\in\mathfrak S_m}\sigma F$,
and $(K_m^0)_n$ is the homogeneous component of degree $n$.
We prove \eqref{SYMMETRIZATION} by a decreasing induction on $n$.
If ${\rm deg}\,g$ is large enough,
we can write $g=\sum_{i=1}^Na_ig_i$ with some $a_i\in S_m^0$
and $g_i\in K_m$. The assertion is evident in this case. Now
take $g\in(K_m^0)_n$ and suppose that
\eqref{SYMMETRIZATION} is true for degree higher than $n$.
Set $\bar g={\rm Sym}\,g$. Then there exist $g_{ij}\in(K_m)_n$ such that
\be
g-\bar g=\sum_{i< j}g_{ij},\quad s_{ij}g_{ij}=-g_{ij},
\en
where $s_{ij}\in\mathfrak S_m$ is the transposition of $i$ and $j$.
We have
\be
{\rm Sym}(gF)=\bar g\cdot{\rm Sym}\,F
+\frac12\sum_{i< j}{\rm Sym}\{g_{ij}(1-s_{ij})F\}.
\en
The first term in the right hand side belongs to $S_m^0W$.
Using \eqref{V1}, \eqref{V2} and that
\be
\pi_0\oplus\pi_4\simeq S^2(\pi_2),\quad
\pi_2\simeq \wedge^2\pi_2
\en
one can write
\be
(1-s_{ij})F=(z_i-z_j)F_{i,j}\quad\hbox{for some $F_{ij}\in V_m$}.
\en
By the induction hypothesis we obtain that
\be
{\rm Sym}\{g_{ij}(z_i-z_j)F_{ij}\}\in S_0W.
\en

Second, we show that
\be
V_m\subset W_m+K_m^0V_m.
\en
In fact, we have seen that for any $F\in V_m$
\be
(1-s_{ij})F\in K_m^0V_m.
\en
Therefore, we have
\be
F={\rm Sym}\,F+(F-{\rm Sym}\,F)\in W_m+K_m^0V_m.
\en
\end{proof}
Define a decreasing filtration on $W_m$ by
\be
J^n(W_m)=\sum_{i\geq n}(S_m)_i\cdot W_m.
\en
\begin{lem}
\label{DECREASING}
The associated graded space ${\rm Gr}^n_J(W_m)=J^n(W_m)/J^{n+1}(W_m)$
is isomorphic to a direct sum of $(\pi_2^{*m})^*:$
\be
{\rm Gr}_J(W_m)=S_m\otimes(\pi_2^{*m})^*.
\en
\end{lem}
For the proof of this lemma we prepare a character identity.
\begin{lem}\label{CHARID}
We have the identity 
\begin{equation}
\label{pi2pi3}
\ch_{q,z} L_{0,1}= \sum_{m\ge 0} \frac{q^{m^2}}{(q)_m}
\ch_{q^{-1},z} \pi_2^{*m}.
\end{equation}
\end{lem}
\begin{proof}
Recall the isomorphism \eqref{Demaz}.
We show that
\begin{equation}
\label{pmn}
q^{N^2} \ch_{q^{-1},z} \pi_1^{*2N}=
\sum_{m=0}^N  q^{m^2} \bin{N}{m} \ch_{q^{-1},z} \pi_2^{*m}.
\end{equation}

Using the relation (see (\ref{S_1}))
\begin{equation}
\label{induct}
\ch_{q,z} \pi_1^{*2k_1} * \pi_2^{*k_2}=
  \ch_{q,z} \pi_1^{*2(k_1-1)} * \pi_2^{*(k_2+1)}+
  q^{2k_1+k_2-1} \ch_{q,z} \pi_1^{*2(k_1-1)}* \pi_2^{*k_2}, 
\end{equation}
we rewrite
\be
q^{N^2} \ch_{q^{-1},z} \pi_1^{*2N}
& =& 
q^{N^2} (\ch_{q^{-1},z} \pi_1^{*2(N-1)}*\pi_2+
q^{-2N+1} \ch_{q^{-1},z} \pi_1^{*2(N-1)})
\\
& =& 
q^{N^2} (\ch_{q^{-1},z} \pi_1^{*2(N-2)}*\pi_2^{*2}
+(q^{-2N+1}+q^{-2N+2}) \ch_{q^{-1},z} \pi_1^{*2(N-2)}*\pi_2
\\
& & \quad+ 
q^{-4N+4}\ch_{q^{-1},z} \pi_1^{*2(N-2)}). 
\en
Repeating this procedure $N$ times we obtain (\ref{pmn}).
Obviously, $$\lim_{N\to\infty} \bin{N}{m}=\frac{1}{(q)_m}.$$ Lemma is proved.
\end{proof}
We fix an isomorphism of $\mathfrak{sl}_2$ modules:
$\mathfrak{sl}_2^*\simeq\pi_2$. Then, we have
\begin{lem}\label{DUALITY}
The mapping $\varphi\mapsto F_\varphi$ is an injection
$Q_m^*\rightarrow z_1\cdots z_mW_m$.
\end{lem}
\begin{proof}
The graded action of $x_n$ ($x\in\mathfrak{sl}_2$) on
${\rm Gr}^G(L_{0,1})=\oplus_{m=0}^\infty Q_m$
is mutually commutative and zero if $n\geq0$. The former property implies
that $F_\varphi(z)$ is symmetric with respect to the action \eqref{ACTION},
and the latter implies that the action
of $x_+(z)$ is equal to that of $x(z)=\sum_{n\in\Z}x_{-n}z^n$
on ${\rm Gr}^G(L_{0,1})$. The integrability of the
representation $L_{0,1}$ implies that $e(z)^2=0$. Therefore, by applying
${\rm ad}\,f_0$, we see that the following  operators act as zero on
${\rm Gr}^G(L_{0,1})$:
\be
e_+(z)^2,\ e_+(z)h_+(z),\ 2e_+(z)f_+(z)-h_+(z)^2,\ h_+(z)f_+(z),\ f_+(z)^2.
\en
All these properties imply \eqref{V1} and \eqref{V2}.
\end{proof}
{\it Proof of Lemmas \ref{CHARACTER} and \ref{DECREASING}}.\quad
Proposition \ref{FUSIONDUAL} and Lemma \ref{NATURAL} imply that
for any symmetric polynomial $g$ we have an isomorphism
\be
gW_m/gS_m^0W_m\simeq (\pi_2^{*m})^*.
\en
Suppose that a vector space $B$ is a subspace of a sum of
vector spaces $A_1+\cdots+A_k$. Then, there is a natural surjection.
\bel
A_1/(A_1\cap B)\oplus\cdots\oplus A_k/(A_k\cap B)\rightarrow
(A_1+\cdots+A_k)/B.
\label{QUOTIENT}
\enl
Take a basis $\{g_1,\ldots,g_k\}$ of $(S_m)_n$. Set
$A_i=g_iW_m$ $(1\leq i\leq k)$ and $B=J^{n+1}(W_m)$.
Then we have $A_1+\cdots+A_k=J^n(W_m)$ and $A_i\cap B=g_iS_m^0W_m$. 
Using \eqref{QUOTIENT} in this setting, we obtain
\be
\ch_qL_{0,1}=\sum_{m=0}^\infty\ch_qQ_m
=\sum_{m=0}^\infty\ch_{q^{-1}}Q_m^*\leq\sum_{m=0}^\infty q^m\ch_{q^{-1}}W_m
\leq\sum_{m=0}^\infty\frac{q^{m^2}}{(q)_m}\ch_{q^{-1}}\pi_2^{*m}.
\en
Here we used \eqref{eq:degree-correspondence-W} to obtain the last inequality. 
Because of Lemma \ref{CHARID}, the left end and the right end are equal.
Lemma \ref{CHARACTER} follows from this. Moreover,
for each $n$, we have the isomorphism
\be
\oplus_{i=1}^kg_iW_m/g_iS_m^0W_m\simeq
\oplus_{i=1}^kg_i\otimes(\pi_2^{*m})^*
\simeq {\rm Gr}_J^n(W_m)
\en
in the above notation. Lemma \ref{DECREASING} follows from this.\qed

{\it Proof of Proposition \ref{Qm}}.\quad
Consider the dual $\mathfrak g_-$-action on $Q_m^*$.
The mapping in Lemma \ref{DUALITY} is $\mathfrak g_-$-linear. Set
\bea
Q'_m=U(\mathfrak g_-)(Q_m)^{2m}.
\label{eq:Q-prime-definition}
\ena

Define a coupling between $Q'_m$ and $z_1\cdots z_mW_m$ as follows.
Take $x\cdot w\in Q'_m$ and $g\in z_1\cdots z_mW_m$.
Here $x\in U(\mathfrak g_-)$ and $w\in (Q_m)^{2m}$.
There exists a non-degenerate coupling \cite{FS} between
$(Q_m)^{2m}$ and $z_1\cdots z_m(W_m)^{-2m}$ induced from
the non-degenerate coupling between the subspace of the tensor algebra
over the vector space $\oplus_{i\in\Z_{< 0}}\C e_i$ which is
spanned by $e_{i_1}\cdots e_{i_m}$
$(i_1,\ldots,i_m\in\Z_{< 0})$, and the space of polynomials
in $z_1,\cdots,z_m$ divisible by $z_1\cdots z_m$:
\be
\langle e_{i_1}\cdots e_{i_m},z_1^{n_1}\cdots z_m^{n_m}\rangle
=\prod_{j=1}^m\delta_{i_j+n_j,0}.
\en
Using this coupling we define
\be
\langle x\cdot w,g\rangle=\langle w,S(x)\cdot g\rangle.
\en
Here $S(x)$ is the antipode of $U(\mathfrak g_-)$.
Using the scaling operator we see that the above coupling
further induces a coupling
\bel
{\rm Gr}^J_n(Q'_m)\times{\rm Gr}^n_J(W_m)\rightarrow\C.
\label{COUPLING}
\enl
The $\mathfrak g_-$-module $\pi_2^{*m}$ is cyclic and generated by
the highest weight vector:
\be
\pi_2^{*m}\subset U(\mathfrak g_-)\cdot(\pi_2^{*m})^{2m}.
\en
The dual module $(\pi_2^{*m})^*$ is cocyclic:
\be
U(\mathfrak g_-)\cdot w\supset\bigl((\pi_2^{*m})^*\bigr)^{2m}
\hbox{ for any $0\neq w\in(\pi_2^{*m})^*$.}
\en
Now we know that ${\rm Gr}^n_J(W_m)$ is a direct sum of cocyclic modules 
{}from Lemma \ref{DECREASING}, 
and $Q_{m}'$ is a cyclic module by definition 
\eqref{eq:Q-prime-definition}. 
Therefore, the coupling \eqref{COUPLING} defines an inclusion
\bel
{\rm Gr}^n_J(W_m)\subset{\rm Gr}^J_n(Q'_m)^*.
\label{INCLUSION}
\enl
Comparing this result with Lemmas \ref{DECREASING} and \ref{CHARID},
we obtain the equality
\be
Q'_m=Q_m
\en
and that the inclusion \eqref{INCLUSION} is bijective.\qed

\subsection{Fusion filtration.}
In this section we conclude our discussion on the fusion filtration
on the tensor product of $\widehat{\mathfrak{sl}}_2$-modules.
\begin{lem}\label{lem:Qsum}
The following formula holds\,$:$ 
\bea
& & \ch_q I_{k+1}^l \bigl({\rm Gr}^{J}(Q_m)\otimes\pi_j\bigr)
\label{m}\\ 
& & \quad 
=
\frac{q^{m^2}}{(q)_m}
\sum_{\la\in\Z} q^{-(k+3)\la^2+(l+1)\la}
\Bigl(
\ch_{q^{-1}} (\pi_2^{*m})^{2(k+3)\la-l+j}-\ch_{q^{-1}} (\pi_2^{*m})^{2(k+3)\la-l-j-2}
\Bigr).
\nonumber
\ena
For $j\equiv l\bmod 2$, we have
\bea
& & 
q^{\Delta(j,k)-\Delta(l,k+1)}
\ch_qM^{(k+2,k+3)}_{j+1,l+1}
=\sum_{m\ge 0}\ch_qI_{k+1}^l({\rm Gr}^{J}(Q_m)\otimes\pi_j).
\label{Qsum}
\ena
\end{lem}
\begin{proof}
In view of Proposition \ref{Qm} and 
the vanishing of homology \eqref{vanish3}, 
we can apply Corollary \ref{Euler2} to 
$V={\rm Gr}^{J}(Q_m)\otimes \pi_j$. 
A simple calculation leads to formula \eqref{m}. 
To show \eqref{Qsum}, we apply the identity of characters
Lemma \ref{lem:rocha} taking 
$p=k+2$, $p'=k+3$, $r=b=j+1$, $a=l+1$ 
and noting the relation
\bea
S_{m,\alpha}(q)=\ch_{q^{-1}}(\pi_2^{*m})^\alpha. 
\label{eq:pi2-to-S}
\ena
Comparing the result with \eqref{m} we obtain \eqref{Qsum}. 
\end{proof}

Let $F_m$ be the fusion filtration \eqref{fusfiltr} on the tensor product
$L_{0,1}(1)\T L_{j,k}(0)$.
\begin{thm}\label{mth}
We have 
\be
& & F_m/F_{m-1}\simeq I_{k+1} (Q_m\T\pi_j).
\en
For the character we have 
\bea
\ch_q I_{k+1}^l \bigl(Q_m\otimes\pi_j\bigr)
=\ch_q I_{k+1}^l \bigl({\rm Gr}^{J}(Q_m)\otimes\pi_j\bigr),
\label{eq:character-identity-GrQ}
\ena
which is given explicitly by \eqref{m}. 
\end{thm}
\begin{proof}
Since there is a surjection
$I_{k+1}(Q_m\otimes\pi_j)\to F_m/F_{m-1}$, 
we are to show that both sides have the same character. 
We note that
\begin{multline}
\label{le}
\ch_q I_{k+1}^l (L_{0,1}\otimes\pi_j)\le \sum_{m\ge 0}
  \ch_q I_{k+1}^l (Q_m\otimes\pi_j)\le
  \sum_{m\ge 0} \ch_q I_{k+1}^l ({\rm Gr}^{J}(Q_m)\otimes\pi_j).
\end{multline}
On the other hand, 
$I_{k+1}(L_{0,1}\otimes\pi_j)=L_{0,1}\T L_{j,k}$ 
implies that 
$$
q^{\triangle(j,k)-\triangle(l,k+1)}
\ch_q I_{k+1}^l (L_{0,1}\T\pi_j)
=\ch_q M_{j+1,l+1}^{k+2,k+3} 
$$
(the factor $q^{\triangle(j,k)-\triangle(l,k+1)}$ comes from the
difference of $L_0$- and $d$-gradings).
Because of Lemma \ref{lem:Qsum},  
we obtain equalities in (\ref{le}).
\end{proof}

Finally, let us discuss the fusion product 
\be
L_{1,1}*L_{j,k}=
\bigoplus_{0\le l\le k+1\atop l\not\equiv j\bmod2} 
{\rm Gr}^E(M^{(k+2,k+3)}_{j+1,l+1})\T L_{l,k+1}.
\en 
Let $\tilde{F}_m$ be the fusion filtration \eqref{fusfiltr}
on the tensor product $L_{1,1}(1) \T L_{j,k}(0)$
of $\slth\T\C[u]$-modules,  
and set 
$$
\tilde{F}_m/\tilde{F}_{m-1}=\bigoplus_{0\le l\le k+1\atop l\not\equiv j\bmod2} 
\tilde{M}^{(k+2,k+3)}_{j+1,l+1,m}\T L_{l,k+1}.
$$

\begin{prop}\label{prop:i=0-to-i=1}
$\ch_q \tilde{M}^{(k+2,k+3)}_{j+1,l+1,m}
=\ch_q{\rm Gr}^E_m(M^{(k+2,k+3)}_{k-j+1,k-l+2}).$
\end{prop}
\begin{proof}
Consider an automorphism 
$\imath:\slth\to \slth$ defined by
$E_i\mapsto E_{1-i}$, $F_i\mapsto F_{1-i}$, 
where $E_0, E_1,F_0, F_1$ are
the Chevalley generators of $\slth$. 
Then we have an isomorphism of
$\slth$-modules
$L_{j,k}\to \tilde L_{k-j,k}$, where 
the action of $\slth$ on the right hand
side is a composition of $\imath$ and 
of the natural action. This proves
our proposition. 
\end{proof}

Let us summarize the conclusion. 
\begin{thm}\label{MAIN}
For the unitary series $M^{(k+2,k+3)}_{r,s}$,  
the character of the (1,3) filtration \eqref{13filt} 
is given by the formula introduced in Lemma \ref{lem:rocha}:
\bea
\ch_q{\rm Gr}^E_m(M^{(k+2,k+3)}_{r,s})=
\frac{q^{\triangle_{r,s}}}{(q)_m} 
I^{(k+2,k+3)}_{r,s,r+i,m}(q).
\label{eq:the-last-formula}
\ena
Here $i=0,1$ is given by $i\equiv r-s\bmod 2$. 
\end{thm}

\begin{proof}
First we give the proof in the case of $i=0$. 
{}From Proposition \ref{fusphi} we have 
\begin{eqnarray*}
\bigoplus_{1\le s\le k+2\atop s\equiv r \bmod 2}
{\rm Gr}^E_m(M^{(k+2,k+3)}_{r,s})\T L_{s-1,k+1}=
(L_{0,1}* L_{r-1,k})_m=F_{m}/F_{m-1}. 
\end{eqnarray*}
Here $F_{m}$ is the fusion filtration on the tensor product 
$L_{0,1}(0) \otimes L_{r-1, k}(1)$. 
{}From Theorem \ref{mth} we have  
\begin{eqnarray*}
{\rm Gr}^E_m(M^{(k+2,k+3)}_{r,s})=
\mathrm{Hom}_{\g'}(L_{s-1,k}, F_{m}/F_{m-1})
\simeq I_{k+1}^{s-1}(Q_{m} \otimes \pi_{r-1}). 
\end{eqnarray*}
Take the character of the both ends above 
and use \eqref{eq:character-identity-GrQ}. 
Then we obtain 
\begin{eqnarray*}
{\rm ch}_{q}{\rm Gr}^E_m(M^{(k+2,k+3)}_{r,s})=
{\rm ch}_{q}I_{k+1}^{s-1}({\rm Gr}^{J}(Q_{m}) \otimes \pi_{r-1}). 
\end{eqnarray*}
The character in the right hand side is given by \eqref{m}.
Then, from the relation \eqref{eq:pi2-to-S}, 
we finally get the equality \eqref{eq:the-last-formula}. 
To show the theorem in the case of $i=1$, use the automorpism of $\widehat{\mathfrak{sl}}_{2}$
given in the proof of Proposition \ref{prop:i=0-to-i=1}.
\end{proof}

\bigskip 
\noindent
{\it Acknowledgments.}\quad
Research of BF is partially supported by RFBR Grants 04-01-00303 and
05-01-01007, INTAS 03-51-3350, NSh-2044.2003.2 and RFBR-JSPS Grant 05-01-02934YaFa.
Research of EF is partially supported by the RFBR Grant 03-01-00167 and LSS 4401.2006.2.
Research of MJ is supported by 
the Grant-in-Aid for Scientific Research B2--16340033. 
Research of TM is supported by 
the Grant-in-Aid for Scientific Research A1--13304010. 
Research of YT is supported by 
the Grant-in-Aid for Young Scientists B--17740089.

\newcounter{a}
\setcounter{a}{2}

\end{document}